\documentclass[final,hidelinks,onefignum,onetabnum]{siamart251216}

\usepackage{amssymb,mathtools}
\usepackage{booktabs,array}
\usepackage{threeparttable}
\usepackage{enumitem,xcolor}
\let\SIAMthebibliography\thebibliography
\let\SIAMendthebibliography\endthebibliography
\usepackage[numbers]{natbib}
\let\thebibliography\SIAMthebibliography
\let\endthebibliography\SIAMendthebibliography

\newsiamthm{assumption}{Assumption}
\newsiamremark{remark}{Remark}
\theoremstyle{plain}
\theoremheaderfont{\normalfont\sc}
\theorembodyfont{\normalfont}
\theoremseparator{.}
\theoremsymbol{}
\newtheorem{algo}[theorem]{Algorithm}

\DeclareMathOperator*{\argmin}{arg\,min}
\DeclareMathOperator*{\argmax}{arg\,max}
\DeclareMathOperator{\dist}{dist}
\DeclareMathOperator{\ri}{ri}
\newcommand{\R}{\mathbb{R}}
\newcommand{\E}{\mathsf{E}}
\newcommand{\ip}[2]{\langle #1,#2\rangle}
\newcommand{\norm}[1]{\lVert #1\rVert}
\newcommand{\dnorm}[1]{\lVert #1\rVert_{*}}
\newcommand{\Xs}{X_\star}
\newcommand{\fs}{f_\star}
\newcommand{\xs}{x_\star}
\newcommand{\NC}{N_C}
\newcommand{\lone}{\ell_1}
\newcommand{\linf}{\ell_\infty}

\headers{Broximal gradient descent in general norms}{Peter Richt\'arik and Xun Qian}

\title{Broximal Gradient Descent in General Norms: Convergence without Fej\'er Monotonicity}

\author{Peter Richt\'arik\thanks{King Abdullah University of Science and Technology (KAUST), Thuwal, Saudi Arabia (\mbox{\email{peter.richtarik@kaust.edu.sa}}, \mbox{\email{xun.qian@kaust.edu.sa}}).}
\and Xun Qian \footnotemark[1]}

\ifpdf
\hypersetup{
  pdftitle={Broximal Gradient Descent in General Norms: Convergence without Fejer Monotonicity},
  pdfauthor={Peter Richtarik and Xun Qian}
}
\fi

\begin{document}

\maketitle

\begin{abstract}
We study Broximal Gradient Descent in a general norm for smooth optimization over a closed convex set. Each update minimizes the linear model over the intersection of the feasible set and a norm ball. Unlike the Euclidean case with suitably restricted radii, the distance to a solution can increase, even when the squared norm is strongly convex. We instead use a comparator inequality and an observable stationarity certificate. For convex objectives, prescribed radii yield value convergence, including escape when the infimum is unattained, and a last-iterate $O(K^{-1/3})$ bound without a distance-bound assumption. A level-set distance bound gives an all-iteration $O(1/K)$ rate with calibrated radii and an asymptotic $O(\log K/K)$ rate with radii using no problem constant. Strong convexity and power growth give explicit faster envelopes and convergence of the iterates. Without convexity, a constrained Polyak--{\L}ojasiewicz inequality gives linear convergence with prescribed radii. An observable radius-halving rule needs a valid smoothness constant but neither the PL modulus nor an initial gap bound; warm starting gives logarithmic total oracle complexity in inverse accuracy. The analysis separates established local-oracle ingredients from the control of level sets, radius calibration, and accuracy-based oracle counts.
\end{abstract}



\section{Introduction}\label{sec:intro}

We study the constrained problem
\[
\min_{x\in C}f(x),
\]
where $C$ is a nonempty closed convex subset of a finite-dimensional real normed space $(\E,\norm\cdot)$ and $f$ is differentiable, by a method that respects both the constraint and the geometry of the norm: at the current point $x_k$ it replaces $f$ by its linear model and minimizes that model over the part of $C$ that lies within distance $r_k$ of $x_k$.

\begin{algo}[BroxGD]\label{alg:lmo}
Given $x_0\in C$ and radii $r_k>0$, for $k=0,1,2,\dots$ compute
\begin{equation}\label{eq:update}
x_{k+1}\in\argmin_{u\in C:\ \norm{u-x_k}\le r_k}\ \ip{\nabla f(x_k)}{u-x_k}.
\end{equation}
\end{algo}

BroxGD is the Euclidean method of Richt{\'a}rik et~al.~\cite{richtarik2026broximal} extended to a chosen norm, and the linearized ball-proximal update of Gruntkowska and Richt{\'a}rik~\cite{gruntkowska2025non}. Its radius is the only step parameter. It avoids projection, but requires a local linear solve: an LP when $C$ and the norm ball are polyhedral. For $C=\E$, it is normalized steepest descent: normalized gradient descent in Euclidean geometry, sign descent in $\linf$, and the orthogonalized matrix update in spectral-norm geometry. In $\lone$, one output is $x_k-r_k\operatorname{sign}(\nabla_{i_k}f(x_k))e_{i_k}$ with $i_k\in\argmax_i|\nabla_if(x_k)|$, a Gauss--Southwell coordinate step of prescribed length \cite{nutini2015coordinate}. Ties allow other outputs, and coupled constraints can require several coordinates to move.

\paragraph{Applications} Closed convex constraints occur in statistical estimation, signal processing, and learning: examples include constrained likelihood estimation \cite{hathaway1985constrained}, real-time convex optimization \cite{mattingley2010real}, and fairness-constrained training \cite{donini2018empirical,cotter2019optimization}. Norm-constrained training also motivates choosing different geometries for different parameter groups \cite{pethick2025training}. These examples motivate the oracle, not a general claim of tractability: nonlinear training constraints need not be convex, and even a convex local subproblem can be expensive. Our guarantees concern closed convex $C$ and count oracle calls.

\paragraph{The question} For admissible Euclidean radii, the decrease $\norm{x_{k+1}-\xs}_2^2\le\norm{x_k-\xs}_2^2-r_k^2$ gives an $O(1/\varepsilon)$ best-iterate complexity and linear rates under strong convexity \cite{richtarik2026broximal}. In a general norm, a step can increase both the chosen-norm and Euclidean distance, even when the squared chosen norm is strongly convex (Section~\ref{sec:onestep}). We therefore use two other facts: a linear-model comparison, which convexity turns into objective control, and a computable stationarity certificate from the subproblem's optimality conditions.

\paragraph{Contributions} The local oracle and the residual inequality used below are established ingredients. Our focus is the failure of distance decrease, convergence without attainment, level-set control, prescribed radii with finite-length guarantees, and adaptive oracle complexity. The calibrated $O(1/K)$ bound is a direct consequence of the comparator recurrence. Table~\ref{tab:glance} lists the guarantees; there $\fs:=\inf_Cf$, $\Delta_k:=f(x_k)-\fs$ for finite $\fs$, $U_k$ denotes a prescribed envelope with $\Delta_k\le U_k$, and $a>0$ is an arbitrary scale. In brief:
\begin{enumerate}[label=(\roman*),leftmargin=2em]
\item \emph{Failure of the distance decrease.} Two examples show that a BroxGD step can increase the distance to the solution set: by the full radius in the $\linf$ norm, where the Euclidean distance increases as well, and for every small radius in the $\ell_{3/2}$ norm, whose square is strongly convex (Section~\ref{sec:onestep}).
\item \emph{Convex objectives, general norm, no boundedness.} With prescribed radii, BroxGD converges in value on every closed convex $C$, including escape to infinity when the infimum is unattained. When a minimizer exists, a fixed horizon gives a last-iterate $O(K^{-1/3})$ bound; a level-set distance bound $D$ gives $\tilde O(LD^2/\varepsilon)$ with a constant radius, an order known for bounded feasible sets \cite{kovalev2025understanding}, and $O(1/K)$ with calibrated decreasing radii. If $D$ bounds the distance on all of $C$, these decreasing radii need only $D$. Without using problem constants, we obtain $O(\log K/K)$ asymptotically.
\item \emph{Regularity.} Strong convexity and power growth of order $p$ give linear, algebraic, or superlinear envelopes for prescribed radii, and under growth the iterates converge to a minimizer along a path of finite length.
\item \emph{Nonconvex objectives.} Without convexity, the $O(K^{-1/2})$ stationarity rate in the dual norm, known from Kovalev~\cite{kovalev2025understanding}, is attained with an observable output rule, and BroxGD converges linearly under a constrained Polyak--{\L}ojasiewicz inequality. A radius-halving test removes the PL constant and the initial gap bound from the algorithm, and a warm start brings its cost to $2K-1$ oracle calls for $K$ updates plus the net number of halvings, which is $O(\log(1/\varepsilon))$ until the gap first falls below $\varepsilon$.
\end{enumerate}

\begin{table}[t]
\caption{The guarantees of this paper, for the update \eqref{eq:update} in a general norm, and the closest earlier results. ``Euclidean'' refers to Richt{\'a}rik et~al.~\cite{richtarik2026broximal}; $D$ is the level-set distance bound of Assumption~\ref{ass:D}, $s$ is the stationarity measure \eqref{eq:normalcone}, and $k_{\rm cap}$ is bounded in Theorem~\ref{thm:warm}.}\label{tab:glance}
\centering\footnotesize
\renewcommand{\arraystretch}{1.0}
\setlength{\tabcolsep}{3.5pt}
\setlength{\aboverulesep}{0.25ex}\setlength{\belowrulesep}{0.45ex}
\begin{threeparttable}
\begin{tabular}{@{}>{\raggedright}p{0.17\textwidth}>{\raggedright}p{0.12\textwidth}>{\raggedright}p{0.275\textwidth}l>{\raggedright\arraybackslash}p{0.25\textwidth}@{}}
\toprule
Setting & Radius $r_k$ & Guarantee & Result & Closest earlier result\\
\midrule
convex & $\sum r_k=\infty$, $\sum r_k^2<\infty$ & $f(x_k)\to\fs$, also when $\fs=-\infty$ & Cor.~\ref{cor:value} & ---\\
\midrule
convex, $\Xs\ne\emptyset$ & $aK^{-2/3}$ & $\Delta_K=O(K^{-1/3})$ & Cor.~\ref{cor:fixed} & Euclidean: $\min_{k\le K}\Delta_k\le\frac{LR^2}{2K}$, radius uses $R\ge\norm{x_0-\xs}_2$\\
\midrule
\quad + distance bound $D$ & $\varepsilon/(LD)$ & $\tilde O(LD^2/\varepsilon)$ calls to $\Delta_K\le\varepsilon$\tnote{b} & Cor.~\ref{cor:constr} & same order when $C$ has diameter $D$ \cite[Cor.~7]{kovalev2025understanding}\\
\midrule
\quad + distance bound $D$ & $2D/(k+\tau)$ & $\Delta_K\le A/(K+\tau)$ for every $K\ge0$; calibration in \eqref{eq:calibrated} & Cor.~\ref{cor:calibrated} & standard scalar recurrence; no monotone descent needed\\
\midrule
\quad + distance bound $D$ & $\frac{a\ln(k+2)}{k+1}$ & $\Delta_K=O(\log K/K)$ for $K\ge\bar k$; $D$ not used & Thm.~\ref{thm:dfree} & ---\\
\midrule
strongly convex & \eqref{eq:scschedule} & $\Delta_k\le(1-\mu/L)^kU_0$ \ ($\mu\le L/2$) & Thm.~\ref{thm:strong} & Euclidean: factor $1-\frac\mu{4L}$\\
\midrule
growth of order $p$ & \eqref{eq:growthrec} & linear ($p=2$), $O(k^{-p/(p-2)})$ ($p>2$), superlinear ($p<2$); $x_k\to\xs$ & Thm.~\ref{thm:growth} & ---\\
\midrule
nonconvex, Lipschitz gradient\tnote{c} & $a/\sqrt K$ & $s(\hat x_K)=O(K^{-1/2})$, observable output $\hat x_K$ & Thm.~\ref{thm:stat} & same rate, no output rule \cite[Thm.~1]{kovalev2025understanding}\\
\midrule
\quad + constrained PL & \eqref{eq:plschedule} & $\Delta_k\le(3L/(3L+\mu))^kU_0$ & Thm.~\ref{thm:pl} & Euclidean:\tnote{a}\ \ factor $1-\frac\mu L$\\
\midrule
\quad + constrained PL & halving \eqref{eq:halving} & linear after finitely many capped steps; accepted-update budget \eqref{eq:halvingbudget} & Thm.~\ref{thm:halving} & Euclidean, strongly convex: backtracking with factor $\xi$, $\Delta_K\le(1+\frac{\xi^2\mu}L)^{-K}\Delta_0$\\
\midrule
\quad + constrained PL & warm-started halving \eqref{eq:warmradius} & $O(\log(1/\varepsilon))$ total oracle calls; explicit dependence in \eqref{eq:warmcalls} & Cor.~\ref{cor:warmcalls} & ---\\
\bottomrule
\end{tabular}
\begin{tablenotes}[flushleft]\footnotesize
\item[a] Under the stronger PL condition on the projected-gradient mapping (Section~\ref{sec:plsub}); the radius $\frac1L\norm{G_{1/L}(x_k)}$ needs a projection.
\item[b] For $0<\varepsilon\le\min\{\Delta_0,LD^2\}$. Settings abbreviate the hypotheses of the cited results. Prescribed schedules use one call per update; the two halving rows distinguish accepted updates from total calls.
\item[c] Gradient Lipschitz continuity is measured in the dual norm, as in Assumption~\ref{as:dual}.
\end{tablenotes}
\end{threeparttable}
\end{table}

\paragraph{Related work} Ferris and Zavriev~\cite{ferris1996linear} analyze \eqref{eq:update} in an arbitrary norm with an Armijo-type radius search and prove linear convergence under the error bound of Luo and Tseng; their test is the ancestor of the observable radius rule of Section~\ref{sec:adaptive}. Garber and Hazan~\cite{garber2016linearly} introduce a relaxed local oracle on polytopes, their local linear optimization oracle, in which the comparison radius and the movement radius differ by a distortion factor, and prove linear convergence under strong convexity. Richt{\'a}rik et~al.~\cite{richtarik2026broximal} give the Euclidean theory of \eqref{eq:update} summarized above. Kovalev~\cite{kovalev2025understanding} studies a non-Euclidean trust-region method with a regularizer, with and without weight decay. Without weight decay and with the indicator of $C$ as regularizer, his update (9) is BroxGD with a constant radius. For this update he proves the $O(K^{-1/2})$ stationarity bound for nonconvex objectives (his Theorem~1; see Lemma~\ref{lem:cert} and Theorem~\ref{thm:stat}) and, for star-convex objectives on a bounded domain of diameter $D$, the complexity $\tilde O(LD^2/\varepsilon)$ (his Corollary~7; see Corollary~\ref{cor:constr}). Pethick et~al.~\cite{pethick2025training} treat norm-constrained linear minimization for training. El~Halabi et~al.~\cite{el2017general} study the generalized proximal gradient step in an arbitrary norm, in which the ball constraint of \eqref{eq:update} is replaced by a squared-norm penalty; here the radius is prescribed directly.

Like the conditional gradient method, BroxGD is projection-free, but its oracle is local. The conditional gradient method has the $O(1/K)$ rate \cite{levitin1966constrained}, in particular with the open-loop step $2/(k+2)$ \cite{dunn1978conditional,jaggi2013revisiting}, and this order is essentially tight for the unmodified method \cite{canon1968tight}; in general, linear rates require further assumptions, such as a solution in the interior of the feasible set, or modifications such as away steps \cite{guelat1986some,lacoste2015global}.

The regularity conditions we use are classical or close to classical: a constrained Polyak--{\L}ojasiewicz inequality, closely related to the proximal Polyak--{\L}ojasiewicz condition \cite{karimi2016linear,bolte2017error}, and quadratic and power growth \cite{drusvyatskiy2018error,necoara2019linear}, with weak sharp minima \cite{burke1993weak} as the case $p=1$. The stationarity measure $s$, the dual-norm distance from $-\nabla f(x)$ to the normal cone of $C$, is the limit as $r\downarrow0$ of the model decrease per unit radius, which at unit radius is the first-order criticality measure of trust-region methods for convex constraints \cite[Chapter~12]{conn2000trust}.

\paragraph{Organization} Section~\ref{sec:onestep} proves the two one-step inequalities and compares them with the Euclidean distance decrease. Section~\ref{sec:convex} treats convex objectives, from no regularity at all to strong convexity and power growth, and Section~\ref{sec:nonconvex} treats stationarity, the constrained Polyak--{\L}ojasiewicz condition, and radius halving without convexity. Section~\ref{sec:related} relates the regularity conditions, and Section~\ref{sec:concl} concludes.

\section{One step: model comparison and certificate at the output}\label{sec:onestep}

\begin{assumption}[geometry and initialization]\label{as:geom}
$\E$ is a finite-dimensional real normed space, $C\subseteq\E$ is nonempty, closed, and convex, and $x_0\in C$. The function $f$ is real-valued and differentiable on an open convex neighborhood of $C$. Derivatives are paired with vectors through $\ip\cdot\cdot$. No smoothness or strict convexity of the norm is assumed.
\end{assumption}

Under Assumption~\ref{as:geom}, the local feasible set in \eqref{eq:update} is nonempty and compact, so a minimizer exists. Each iteration evaluates one gradient and makes one exact local linear-oracle call; all results allow every oracle selection, and oracle counts do not assert that the local subproblem is inexpensive. Throughout, $B(x,r):=\{u:\norm{u-x}\le r\}$, $\dnorm\cdot$ is the dual norm, $\fs:=\inf_Cf$, $\Xs$ is the set of minimizers and $\dist(x,\Xs):=\inf_{z\in\Xs}\norm{x-z}$ when this set is nonempty, and $\Delta(x):=f(x)-\fs$, $\Delta_k:=\Delta(x_k)$ whenever $\fs$ is finite.

\begin{assumption}[convexity]\label{as:conv}
The function $f$ is convex on its open convex neighborhood of $C$.
\end{assumption}

\begin{assumption}[quadratic upper model]\label{as:quad}
For some $L>0$,
\begin{equation}\label{eq:upper}
f(y)\le f(x)+\ip{\nabla f(x)}{y-x}+\frac L2\norm{y-x}^2\qquad(x,y\in C).
\end{equation}
\end{assumption}

Some smoothness of $f$ beyond Lipschitz continuity, such as Assumption~\ref{as:quad}, is needed for the convergence guarantees below. If $f$ is only assumed convex and Lipschitz continuous and the oracle is run with subgradients in place of gradients, BroxGD can fail to converge to $X_\star$ from almost every initial point, whatever the choice of radii; an explicit example in the $\ell_\infty$ norm is given in Proposition~\ref{prop:ce} at the end of this section.

\begin{assumption}[Lipschitz gradient in the dual norm and lower bound]\label{as:dual}
The function $f$ need not be convex. For some $L>0$, its gradient satisfies $\dnorm{\nabla f(y)-\nabla f(x)}\le L\norm{y-x}$ for all $x,y\in C$, where $\dnorm\cdot$ is the dual norm. Also $\fs:=\inf_Cf>-\infty$. Together with Assumption~\ref{as:geom}, this implies the upper model \eqref{eq:upper} by integration along a segment. It is stronger than merely assuming Assumption~\ref{as:quad}.
\end{assumption}

Every result in this paper rests on two inequalities about a single update, one for the comparator route and one for the residual route. We state them first, in a form that does not assume convexity, and show how convexity enters as a corollary. Throughout, $x_{k+1}$ is any output of \eqref{eq:update}, $g_k:=\nabla f(x_k)$, and
\begin{equation}\label{eq:dk}
d_k:=-\ip{g_k}{x_{k+1}-x_k}=-\min_{u\in C,\ \norm{u-x_k}\le r_k}\ip{g_k}{u-x_k}\ \ge0
\end{equation}
is the decrease of the linear model achieved by the oracle; $d_k\ge0$ because $u=x_k$ is feasible.

\begin{lemma}[one-step model comparison]\label{lem:comp}
Let Assumptions~\ref{as:geom} and~\ref{as:quad} hold. Then every output of \eqref{eq:update} satisfies:
\begin{enumerate}[label=(\alph*),leftmargin=1.6em]
\item (linear-model comparison) for every $u\in C$ with $\norm{u-x_k}\le r_k$,
\begin{equation}\label{eq:modelcomp}
f(x_{k+1})\le f(x_k)+\ip{g_k}{u-x_k}+\frac L2r_k^2 ;
\end{equation}
\item (descent) $f(x_{k+1})\le f(x_k)-d_k+\frac L2\norm{x_{k+1}-x_k}^2\le f(x_k)-d_k+\frac L2r_k^2$;
\item (segment comparator) for every $z\in C$ and every $\theta\in[0,1]$ with $\theta\norm{z-x_k}\le r_k$,
\begin{equation}\label{eq:segment}
f(x_{k+1})\le f(x_k)+\theta\,\ip{g_k}{z-x_k}+\frac L2r_k^2 .
\end{equation}
\end{enumerate}
Convexity of $f$ is not assumed. The same inequalities hold for a generic center $x$, radius $r$, and output $q$, with $d:=-\ip{\nabla f(x)}{q-x}$ in place of $d_k$; we use this form whenever a single step is analyzed.
\end{lemma}

\begin{proof}
The upper model and oracle optimality give, for every feasible comparator $u$,
\begin{eqnarray*}
f(x_{k+1})&\overset{\eqref{eq:upper}}\le&f(x_k)+\ip{g_k}{x_{k+1}-x_k}+\frac L2\norm{x_{k+1}-x_k}^2\\
&\overset{\eqref{eq:update}}\le&f(x_k)+\ip{g_k}{u-x_k}+\frac L2r_k^2,
\end{eqnarray*}
which proves (a). Substituting \eqref{eq:dk} into the first line proves (b). For (c), use $u=x_k+\theta(z-x_k)$, which belongs to $C$ and satisfies $\norm{u-x_k}\le r_k$.
\end{proof}

\begin{corollary}[how convexity enters]\label{cor:comp}
Let Assumptions~\ref{as:geom}, \ref{as:conv}, and~\ref{as:quad} hold, and let $\theta\in[0,1]$, $z\in C$ with $\theta\norm{z-x_k}\le r_k$. Then (objective comparison, contraction, and a lower bound on the model decrease)
\begin{align}
f(x_{k+1})&\le f(u)+\omega_k\qquad(u\in C,\ \norm{u-x_k}\le r_k),\qquad \omega_k:=\frac L2r_k^2 ,\label{eq:comp}\\
f(x_{k+1})-f(z)&\le(1-\theta)\bigl(f(x_k)-f(z)\bigr)+\frac L2r_k^2 ,\label{eq:contract}\\
d_k&\ge\theta\bigl(f(x_k)-f(z)\bigr).\label{eq:dlower}
\end{align}
In particular, if $\Xs\ne\emptyset$ and $z$ is a nearest minimizer of $x_k$, then $d_k\ge\min\{1,r_k/\dist(x_k,\Xs)\}\,\Delta_k$, where the minimum is read as $1$ if $x_k\in\Xs$.
\end{corollary}

\begin{proof}
Convexity gives $\ip{g_k}{u-x_k}\le f(u)-f(x_k)$. Inserting this into \eqref{eq:modelcomp} proves \eqref{eq:comp}, and using it with $u=z$ in \eqref{eq:segment} proves \eqref{eq:contract}. For \eqref{eq:dlower}, compare the oracle output with the feasible point $u=x_k+\theta(z-x_k)$:
\[
d_k\overset{\eqref{eq:dk},\eqref{eq:update}}\ge-\theta\ip{g_k}{z-x_k}\ge\theta\bigl(f(x_k)-f(z)\bigr).
\]
Finally choose a nearest minimizer $z$ and $\theta=\min\{1,r_k/\dist(x_k,\Xs)\}$, with the stated zero-distance convention.
\end{proof}

Equations~\eqref{eq:comp}--\eqref{eq:dlower} will control objective gaps. The second route uses the stationarity residual introduced next.

For $y\in C$, define
\begin{equation}\label{eq:normalcone}
\NC(y):=\{n\in\E^*:\ \ip n{z-y}\le0\ \forall z\in C\},\qquad s(y):=\inf_{n\in\NC(y)}\dnorm{\nabla f(y)+n}.
\end{equation}
The normal cone contains the outward normal functionals to $C$. Thus $s(y)$ measures how closely a constraint normal can cancel the derivative; $s(y)=0$ is first-order stationarity.

\paragraph{The model decrease function and $s(x)$} The quantity $d_k$ of Lemma~\ref{lem:comp} is the value at $x=x_k$, $r=r_k$ of the model decrease function
\[
d_x(r):=-\min_{u\in C,\ \norm{u-x}\le r}\ip g{u-x}\qquad(r\ge0),\qquad g=\nabla f(x),
\]
and $s(x)$ is its slope at the origin; we write $d(r)$ for $d_x(r)$ when the center $x$ is fixed, and $d_k:=d_{x_k}(r_k)$ along a run. The next lemma relates $d(r)$ and $s(x)$.  

\begin{lemma}[$d(r)$ and $s(x)$]\label{lem:dr}
Let Assumption~\ref{as:geom} hold and $x\in C$.
(i) $d$ is concave and nondecreasing with $d(0)=0$; hence $d(r)/r$ is nonincreasing in $r$ and $\sup_{r>0}d(r)/r=\lim_{r\downarrow0}d(r)/r$.
(ii) $\displaystyle\sup_{r>0}\frac{d(r)}r=\max_{h\in T_C(x),\ \norm h\le1}\ip{-g}h=\min_{n\in\NC(x)}\dnorm{g+n}=s(x)$, where $T_C(x)=\NC(x)^\circ$ is the tangent cone. In particular $d(r)\le r\,s(x)$ for all $r>0$.
(iii) For every $r>0$: $d(r)=0$ if and only if $s(x)=0$. Consequently either $d\equiv0$ on $[0,\infty)$ (when $s(x)=0$) or $d(r)>0$ for every $r>0$ (when $s(x)>0$).
\end{lemma}
\begin{proof}
\emph{(i)} Let $r_1,r_2\ge0$, $t\in[0,1]$, and let $u_i\in C$ with $\norm{u_i-x}\le r_i$ attain $d(r_i)$. The point $u:=(1-t)u_1+tu_2$ lies in $C$ and satisfies $\norm{u-x}\le(1-t)r_1+tr_2$, so
\[
d\bigl((1-t)r_1+tr_2\bigr)\ \ge\ -\ip g{u-x}=(1-t)\,d(r_1)+t\,d(r_2) ,
\]
and $d$ is concave. It is nondecreasing because the feasible set grows with $r$. Concavity and $d(0)=0$ give $d(tr)\ge t\,d(r)$ for $t\in[0,1]$, i.e.\ $d(r)/r$ is nonincreasing.

\emph{(ii)} Substituting $u=x+rh$,
\[
\frac{d(r)}r=\max\Bigl\{\ip{-g}h:\ \norm h\le1,\ h\in\frac{C-x}r\Bigr\}.
\]
Since $C$ is convex and contains $x$, the sets $(C-x)/r$ increase as $r\downarrow0$, and their union is the cone of feasible directions; by (i), $\sup_{r>0}d(r)/r=\lim_{r\downarrow0}d(r)/r$. The closure of this cone is the tangent cone $T:=T_C(x)$ \cite[Ch.~6]{rockafellar1998variational}; this gives the first equality. For the second, write $N:=\NC(x)$. Conic duality, $N^\circ=T$ \cite[Ch.~6]{rockafellar1998variational}, and Sion's minimax theorem \cite{sion1958general}, applicable because the set $\{\norm h\le1\}$ is compact, give
\begin{eqnarray*}
\min_{n\in N}\dnorm{g+n} &=& \min_{n\in N}\ \max_{\norm h\le1}\ip{-g-n}h\\
&=& \max_{\norm h\le1}\ \inf_{n\in N}\bigl[\ip{-g}h-\ip nh\bigr]\\
&=& \max_{h\in T,\ \norm h\le1}\ip{-g}h ,
\end{eqnarray*}
where the last equality holds because $\inf_{n\in N}(-\ip nh)$ equals $0$ if $h\in N^\circ=T$ and $-\infty$ otherwise. In particular $d(r)\le r\,s(x)$ for every $r>0$.

\emph{(iii)} If $s(x)=0$, then (ii) and monotonicity give $0=d(0)\le d(r)\le r\,s(x)=0$ for every $r>0$. Conversely, let $d(r)=0$ for some $r>0$. For $0<t\le r$, monotonicity gives $0=d(0)\le d(t)\le d(r)=0$, so $d(t)/t=0$ on $(0,r]$, and $s(x)=\lim_{t\downarrow0}d(t)/t=0$ by (i) and (ii). The dichotomy follows: if $s(x)>0$, no $r>0$ can have $d(r)=0$.
\end{proof}

\begin{lemma}[certificate at the output]\label{lem:cert}
Let Assumptions~\ref{as:geom} and~\ref{as:dual} hold, let $q$ be an output of the exact BroxGD oracle \eqref{eq:update} at $x\in C$ with radius $r>0$, and set $d:=-\ip{\nabla f(x)}{q-x}\ge0$. Then
\begin{equation}\label{eq:cert}
s(q)\le Lr+\frac dr ,
\end{equation}
and consequently, with $\Delta(y):=f(y)-\fs$,
\begin{equation}\label{eq:master}
\Delta(q)+r\,s(q)\le\Delta(x)+\frac{3L}2r^2 .
\end{equation}
\end{lemma}
\begin{proof}
\emph{Optimality conditions.} Since $x\in C\cap\operatorname{int}B(x,r)$, the qualification $\ri C\cap\operatorname{int}B(x,r)\ne\emptyset$ holds ($\ri C$ is dense in $C$), and the additivity of subdifferentials \cite[Theorem 23.8]{rockafellar1997convex}, applied to the indicator functions $\iota_C$ and $\iota_{B(x,r)}$, gives
\[
N_{C\cap B(x,r)}(q)=\NC(q)+N_{B(x,r)}(q)\qquad(q\in C\cap B(x,r)) .
\]
By the normal cone of a level set \cite[Corollary~23.7.1]{rockafellar1997convex}, applied to $\varphi(u)=\norm{u-x}$ with Slater point $u=x$,
\[
N_{B(x,r)}(q)=\begin{cases}\{\lambda w:\ \lambda\ge0,\ w\in\partial\norm\cdot(q-x)\}, & \norm{q-x}=r,\\ \{0\}, & \norm{q-x}<r,\end{cases}
\]
In the active-ball case $\norm{q-x}=r>0$, every $w\in\partial\norm\cdot(q-x)$ satisfies $\dnorm w=1$ and $\ip w{q-x}=r$. At $q=x$, the subdifferential is instead the closed dual unit ball; its elements need not have dual norm one. Hence the output $q$ satisfies
\[
\nabla f(x)+\lambda v+n=0,\qquad n\in\NC(q),\ \lambda\ge0,
\]
where $\dnorm v=1$ and $\ip v{q-x}=r$ if $\lambda>0$ (and $\lambda=0$ if the ball is inactive).

\emph{Proof of \eqref{eq:cert}.} Since $x\in C$ and $n\in\NC(q)$, $\ip n{q-x}\ge0$, so
\[
d=-\ip{\nabla f(x)}{q-x}=\lambda\ip v{q-x}+\ip n{q-x}\ \ge\ \lambda r .
\]
Therefore, by Assumption~\ref{as:dual},
\begin{eqnarray*}
s(q) &\overset{\eqref{eq:normalcone}}\le& \dnorm{\nabla f(q)+n}\ =\ \dnorm{\nabla f(q)-\nabla f(x)-\lambda v}\\
&\le& \dnorm{\nabla f(q)-\nabla f(x)}+\lambda\ \le\ Lr+\frac dr .
\end{eqnarray*}

\emph{Proof of \eqref{eq:master}.} Lemma~\ref{lem:comp}(b) gives $\Delta(q)\le\Delta(x)-d+\frac L2r^2$, and \eqref{eq:cert} gives $d\ge r\,s(q)-Lr^2$; combine the two.
\end{proof}

Inequality~\eqref{eq:master} is the residual counterpart of \eqref{eq:contract}. For $C=\E$, $d=r\dnorm{\nabla f(x)}$ and $s=\dnorm{\nabla f}$; Lemma~\ref{lem:comp}(b) then recovers the normalized-steepest-descent inequality $f(q)\le f(x)-r\dnorm{\nabla f(x)}+Lr^2/2$.

Inequality~\eqref{eq:master} is established in the proof of Kovalev's Theorem~1 \cite{kovalev2025understanding}. His regularized local update gives $F(q)+r\dnorm{\nabla f(q)+\hat\nabla R(q)}\le F(x)+3Lr^2/2$, where $F=f+R$ and the subproblem supplies $\hat\nabla R(q)\in\partial R(q)$. Taking $R=\iota_C$ yields \eqref{eq:master}, with the same constant and normal-cone argument. The intermediate certificate \eqref{eq:cert}, not stated there, uses only the observed decrease and radius; it supplies our output selection and radius test. Theorem~\ref{thm:pl} uses the established inequality \eqref{eq:master}; its contribution is the prescribed-radius construction and guarantee.

\paragraph{Comparison with the Euclidean one-step geometry} In the Euclidean norm, Richt{\'a}rik et~al.~\cite[Theorem~3.1]{richtarik2026broximal} prove a stronger one-step statement, which uses neither smoothness nor convexity: if $\xs\in\Xs$ and the radius is admissible, $0<r\le\ip{\nabla f(x)}{x-\xs}/\norm{\nabla f(x)}_2$, then every output satisfies
\begin{equation}\label{eq:fejer}
\norm{q-\xs}_2^2\ \le\ \norm{x-\xs}_2^2-r^2 ,
\end{equation}
from which their Euclidean rates follow. In a general norm the analogue of \eqref{eq:fejer} fails even in its weakest form. Take $C=\E=\R^2$ with the $\linf$ norm, $f(x)=\frac12x^\top\bigl(\begin{smallmatrix}1&-2\\-2&8\end{smallmatrix}\bigr)x$, so that $\Xs=\{0\}$, and the center $x=(1,1)$. Then $\nabla f(x)=(-1,6)$, the natural admissibility condition $r\le\ip{\nabla f(x)}{x-\xs}/\dnorm{\nabla f(x)}=\frac57$ holds for $r=\frac12$, and the unique output is $q=x-r\operatorname{sign}(\nabla f(x))=(\frac32,\frac12)$. The objective drops from $f(x)=\frac52$ to $f(q)=\frac58$, yet
\[
\norm{q-\xs}_\infty=\tfrac32=\norm{x-\xs}_\infty+r :
\]
the distance to the solution in the chosen norm increases by the full radius instead of decreasing. The failure is not caused by the nonsmoothness of the $\linf$ norm, and it persists as $r\downarrow0$. In the $\ell_{3/2}$ norm, whose square is strongly convex, take $f(x)=\frac12x^\top\bigl(\begin{smallmatrix}9&20\\20&48\end{smallmatrix}\bigr)x$, so that $\Xs=\{0\}$, and $x=(4,-1)$. Then $\nabla f(x)=(16,32)$ and the output is $q=x-r\,9^{-2/3}(1,4)$, so that
\[
\frac{d}{dr}\,\norm{q-\xs}_{3/2}\Big|_{r=0}=\frac29>0 :
\]
for every small $r>0$ the distance to the solution increases, although $f(q)<f(x)$. Thus the distance to the solution set in the chosen norm is not monotone along BroxGD steps. Neither is the Euclidean distance: in the $\linf$ example above, $\norm{q-\xs}_2^2=\frac52>2=\norm{x-\xs}_2^2$. These examples rule out the Fej\'er-type argument; they do not show that no other Lyapunov function exists. This is why the analysis below compares objective values through \eqref{eq:segment} and controls the distance to the solution set by separate means; how that distance is controlled is what separates the general-norm rates of Section~\ref{sec:convex} from their Euclidean counterparts.

\paragraph{Smoothness beyond Lipschitz continuity cannot be dropped} We return to the remark after Assumption~\ref{as:quad}. The following example, in which the $\linf$ oracle keeps only the signs of the subgradient, is a variant of Counterexample~2 of Karimireddy et~al.~\cite{karimireddy2019error} for sign descent; Parshakova et~al.~\cite{parshakova2026muon} and Lai et~al.~\cite{lai2026nonsmooth} lift the same mechanism to matrices.

Let $\E=\R^2$ with the $\linf$ norm (dual norm $\lone$), $C=\E$, and fix $\varepsilon\in(0,1)$. Consider
\begin{equation}\label{eq:ce-f}
	f(x)=|x_1+x_2|+\varepsilon\bigl(|x_1|+|x_2|\bigr).
\end{equation}
Then $f$ is convex and piecewise linear, $X_\star=\{0\}$, $f_\star=0$, and every subgradient satisfies $\|g\|_1\le2+2\varepsilon$, so $f$ is $(2+2\varepsilon)$-Lipschitz in the $\ell_\infty$ norm. Since
\[
\partial f(x)=\operatorname{SGN}(x_1+x_2)\,(1,1)+\varepsilon\bigl(\operatorname{SGN}(x_1),\operatorname{SGN}(x_2)\bigr),
\qquad
\operatorname{SGN}(t)=\begin{cases}\{\operatorname{sgn}t\},&t\ne0,\\ [-1,1],&t=0,\end{cases}
\]
the map
\begin{equation}\label{eq:ce-sel}
	g(x):=\bigl(\operatorname{sgn}(x_1+x_2)+\varepsilon\operatorname{sgn}(x_1),\ \operatorname{sgn}(x_1+x_2)+\varepsilon\operatorname{sgn}(x_2)\bigr),\qquad\operatorname{sgn}(0):=1,
\end{equation}
is a subgradient selection: $g(x)\in\partial f(x)$ for every $x$, and $g(x)=\nabla f(x)$ wherever $f$ is differentiable. We run BroxGD with $g(x_k)$ in place of $\nabla f(x_k)$.

\begin{proposition}[a convex Lipschitz function on which BroxGD does not converge]\label{prop:ce}
	Let $f$ be as in \eqref{eq:ce-f}, let $g$ be as in \eqref{eq:ce-sel}, and let
	\[
	x_{k+1}\in\arg\min\bigl\{\langle g(x_k),u-x_k\rangle:\ \|u-x_k\|_\infty\le r_k\bigr\}.
	\]
	Then for every $x_0\in\mathbb{R}^2$ and every sequence of radii $r_k>0$, which may depend on $x_0$, on $k$, and on all previous iterates and subgradients,
	\[
	x_{k,1}-x_{k,2}=x_{0,1}-x_{0,2}
	\qquad\text{and}\qquad
	f(x_k)-f_\star\ \ge\ \varepsilon\,|x_{0,1}-x_{0,2}|
	\qquad\text{for all }k\ge0.
	\]
	In particular, if $x_{0,1}\ne x_{0,2}$, then $f(x_k)\not\to f_\star$ and $\operatorname{dist}(x_k,X_\star)\ge|x_{0,1}-x_{0,2}|/2$ for every choice of radii.
\end{proposition}
\begin{proof}
	Write $\sigma_k:=\operatorname{sgn}(x_{k,1}+x_{k,2})\in\{\pm1\}$. By \eqref{eq:ce-sel}, $g(x_k)=(\sigma_k+\varepsilon a_k,\ \sigma_k+\varepsilon b_k)$ with $a_k,b_k\in\{\pm1\}$. Since $\varepsilon<1$, both components are nonzero and have the sign $\sigma_k$, so the linear minimization over the $\ell_\infty$ ball has the unique solution
	\[
	x_{k+1}=x_k-r_k\operatorname{sgn}\bigl(g(x_k)\bigr)=x_k-\sigma_kr_k(1,1),
	\]
	where $\operatorname{sgn}$ acts componentwise. Both coordinates change by the same amount, hence $x_{k+1,1}-x_{k+1,2}=x_{k,1}-x_{k,2}$, and by induction $x_{k,1}-x_{k,2}=x_{0,1}-x_{0,2}$ for all $k$. Finally, $f(x)\ge\varepsilon(|x_1|+|x_2|)\ge\varepsilon|x_1-x_2|$ and $\|x\|_\infty\ge|x_1-x_2|/2$.
\end{proof}

\begin{remark}\label{rem:ce}
Proposition~\ref{prop:ce} fixes the selection \eqref{eq:ce-sel}. For almost every translate of $f$, its conclusion holds for every subgradient selection. Precisely, consider a deterministic radius rule that uses $x_0$, $k$, its own parameters and the observed subgradients but not function values; this covers the prescribed schedules of Section~\ref{sec:convex} and the halving test of Section~\ref{sec:adaptive}, whose test value $d(r)=r\|g\|_1$ depends on the center only through the observed subgradient $g$. Fix such a rule and $x_0$, and for $c\in\R^2$ let $f_c:=f(\cdot-c)$. Then there is a set $N\subset\R^2$ of Lebesgue measure zero such that, for every $c\notin N$, the iterates generated on $f_c$ are the same for every subgradient selection, $x_{k,1}-x_{k,2}=x_{0,1}-x_{0,2}$ for all $k$, and
\[
f_c(x_k)-\min f_c\ \ge\ \varepsilon\,\bigl|(x_{0,1}-c_1)-(x_{0,2}-c_2)\bigr|\ >\ 0\qquad\text{for all }k\ge0 .
\]
The proof is given in Appendix~\ref{app:ce}. For rules that evaluate $f$, the dependence on the selection cannot be removed on this example: an exact line search along the step direction lands on the line $x_1+x_2=0$, after which a selection with mixed signs makes progress. In the Euclidean norm this obstruction is absent for admissible radii: the Fej\'er-type inequality \eqref{eq:fejer} of Richt{\'a}rik et~al.~\cite[Theorem~3.1]{richtarik2026broximal} holds with a subgradient in place of the gradient.
\end{remark}

\section{Convex objectives}\label{sec:convex}

We proceed from no assumption beyond convexity and the quadratic model (Section~\ref{sec:genconvex}) to strong convexity and power growth, each of which strengthens the control of $\norm{x_k-\xs}$ and with it the rate.

\subsection{General convex objectives}\label{sec:genconvex}

\subsubsection{A comparator theorem without attainment}\label{sec:comparator}

Fix any $z\in C$ and write
\begin{equation}\label{eq:Bz}
R^z:=\norm{x_0-z},\qquad B_k(z):=R^z+\sum_{j<k}r_j,\qquad \Delta_k(z):=f(x_k)-f(z).
\end{equation}
We choose the comparison point that makes these weighted differences telescope.

\begin{theorem}[signed comparator bound]\label{thm:comparator}
Under Assumptions~\ref{as:geom} and~\ref{as:conv}, let $x_{k+1}\in C$ have step length at most $r_k$ and satisfy the comparison inequality in \eqref{eq:comp} with any nonnegative errors $\omega_k$. For every $z\in C$ and $K\ge1$,
\begin{equation}\label{eq:signed}
B_K(z)\Delta_K(z)\le R^z\Delta_0(z)+\sum_{k=0}^{K-1}B_{k+1}(z)\omega_k .
\end{equation}
In particular, BroxGD under Assumption~\ref{as:quad} obeys this bound with $\omega_k=Lr_k^2/2$.
\end{theorem}
\begin{proof}
The triangle inequality gives $\norm{x_k-z}\le R^z+\sum_{j<k}r_j=B_k(z)$. Set $\theta_k:=r_k/B_{k+1}(z)\in(0,1]$ and $u_k:=(1-\theta_k)x_k+\theta_kz\in C$. Then $\norm{u_k-x_k}=\theta_k\norm{z-x_k}\le\theta_kB_k(z)\le r_k$, so the comparison inequality at $u_k$ and convexity of $f$ give, as in \eqref{eq:contract},
\[
\Delta_{k+1}(z)\overset{\eqref{eq:comp}}\le(1-\theta_k)\Delta_k(z)+\omega_k\overset{\eqref{eq:Bz}}=\frac{B_k(z)}{B_{k+1}(z)}\Delta_k(z)+\omega_k ,
\]
where $1-\theta_k=B_k(z)/B_{k+1}(z)$ because $B_{k+1}(z)=B_k(z)+r_k$. Multiplying by $B_{k+1}(z)>0$,
\[
B_{k+1}(z)\Delta_{k+1}(z)\le B_k(z)\Delta_k(z)+B_{k+1}(z)\,\omega_k ,
\]
and summing over $k<K$ with $B_0(z)=R^z$ gives \eqref{eq:signed}. For BroxGD, Corollary~\ref{cor:comp} provides the comparison inequality with $\omega_k=Lr_k^2/2$.
\end{proof}

\begin{corollary}[value convergence and escape]\label{cor:value}
Under Assumptions~\ref{as:geom}, \ref{as:conv}, and~\ref{as:quad}, suppose $\sum_{k\ge0}r_k=\infty$ and $\sum_{k\ge0}r_k^2<\infty$. Then
\begin{equation}\label{eq:valueconv}
f(x_k)\longrightarrow\fs:=\inf_{x\in C}f(x),
\end{equation}
including $\fs=-\infty$. If the infimum is unattained, then $\norm{x_k}\to\infty$.
\end{corollary}
\begin{proof}
Fix $z\in C$ and $M<K$. Dividing \eqref{eq:signed} by $B_K(z)$ and using $B_{k+1}(z)\le B_K(z)$ for $k\ge M$,
\[
\Delta_K(z)\ \overset{\eqref{eq:signed}}\le\ \frac{R^z\Delta_0(z)+\sum_{k<M}B_{k+1}(z)\,\omega_k}{B_K(z)}+\frac L2\sum_{k\ge M}r_k^2 .
\]
Since $B_K(z)\to\infty$, letting $K\to\infty$ and then $M\to\infty$ gives $\limsup_Kf(x_K)\le f(z)$. Taking the infimum over $z\in C$ proves \eqref{eq:valueconv}; when $\fs=-\infty$, choose $z$ below any prescribed value. If the infimum is unattained and a subsequence of $(x_k)$ were bounded, it would have a cluster point in $C$, which by continuity of $f$ would be a minimizer (or contradict $f(x_k)\to-\infty$). Hence $\norm{x_k}\to\infty$.
\end{proof}

The practical schedule $r_k=a/(k+1)^\beta$, with $a>0$ and $1/2<\beta\le1$, needs neither a solution nor an initial gap bound. The conclusion concerns the last objective value; it does not assert that these values decrease monotonically.

\subsubsection{Rates when a minimizer exists}\label{sec:rates}

To express a finite-horizon rate using the initial distance, we now choose an optimal comparator. 

\begin{assumption}[attainment]\label{as:att}
The set $\Xs:=\argmin_{x\in C}f(x)$ is nonempty. Fix $\xs\in\Xs$ when a single optimal comparator is needed, and define
\begin{equation}\label{eq:att}
\Delta_k:=f(x_k)-\fs,\qquad R_0:=\norm{x_0-\xs},\qquad B_k:=R_0+\sum_{j<k}r_j .
\end{equation}
\end{assumption}

Under Assumptions~\ref{as:geom}, \ref{as:conv}, \ref{as:quad}, and~\ref{as:att}, Theorem~\ref{thm:comparator} gives
\begin{equation}\label{eq:attbound}
\Delta_K\overset{\eqref{eq:signed}}\le\frac{R_0\Delta_0+\frac L2\sum_{k<K}B_{k+1}r_k^2}{B_K}.
\end{equation}

\begin{corollary}[a fixed budget of calls]\label{cor:fixed}
Under Assumptions~\ref{as:geom}, \ref{as:conv}, \ref{as:quad}, and~\ref{as:att}, fix $K\ge1$ and $a>0$. With $r_k=aK^{-2/3}$ for $k<K$,
\begin{equation}\label{eq:fixedrate}
\Delta_K\le\frac{R_0\Delta_0}{R_0+aK^{1/3}}+\frac{La^2}{2K^{1/3}}\le\Bigl(\frac{R_0\Delta_0}a+\frac{La^2}2\Bigr)K^{-1/3}.
\end{equation}
\end{corollary}

\begin{proof}
The constant radius gives $B_K=R_0+aK^{1/3}$ and $\sum_{k<K}r_k^2=a^2K^{-1/3}$. Since $B_{k+1}\le B_K$,
\[
\Delta_K\overset{\eqref{eq:attbound}}\le\frac{R_0\Delta_0}{B_K}+\frac L2\sum_{k<K}r_k^2=\frac{R_0\Delta_0}{R_0+aK^{1/3}}+\frac{La^2}{2K^{1/3}}.
\]
Use $B_K\ge aK^{1/3}$ for the second inequality.
\end{proof}

For constant $r$, the simpler bound is $R_0\Delta_0/(Kr)+LKr^2/2$. Balancing accumulated progress against model error explains $r\propto K^{-2/3}$. The parameter $a$ is arbitrary; optimizing its constant requires $R_0\Delta_0$.

\paragraph{An uninterrupted anytime schedule} Under Assumptions~\ref{as:geom}, \ref{as:conv}, \ref{as:quad}, and~\ref{as:att}, choose $r_k=a/(k+1)^{2/3}$. Then
\begin{equation}\label{eq:anytime}
\Delta_K\le\frac{R_0\Delta_0+2LR_0a^2+\frac32La^3(1+\log K)}{R_0+3a((K+1)^{1/3}-1)}=O\Bigl(\frac{\log K}{K^{1/3}}\Bigr).
\end{equation}
Indeed, with $S_j:=\sum_{k<j}r_k$, integral comparison gives $3a((j+1)^{1/3}-1)\le S_j\le3aj^{1/3}$; moreover $\sum_{j=1}^Kj^{-4/3}\le4$ and $\sum_{j=1}^Kj^{-1}\le1+\log K$, and substituting these bounds into \eqref{eq:attbound} gives \eqref{eq:anytime}. These upper bounds do not rule out an $O(1/K)$ arbitrary-norm result; the Euclidean squared-distance argument behind the $O(1/\varepsilon)$ target-radius bound is not available in a general norm (Section~\ref{sec:onestep}).

\subsubsection{Last-iterate \texorpdfstring{$O(1/K)$}{O(1/K)} rates under a level-set distance bound}\label{sec:levelset}

The fixed-comparator bound of Corollary~\ref{cor:fixed} lets the distance estimate grow with the accumulated radii. A uniform distance bound instead gives the recurrence $\Delta_{k+1}\le(1-r_k/D)\Delta_k+Lr_k^2/2$. The issue is keeping the iterates in the region where that bound holds, not requiring their objective values to decrease. Lemma~\ref{lem:inv} proves this invariance for sufficiently small radii. Constant radii then give a target-accuracy guarantee, while calibrated decreasing radii give $O(1/K)$ for every iteration. Finally, Theorem~\ref{thm:dfree} removes problem constants from the radii, with a logarithmic factor and a potentially large onset.

\begin{assumption}[bounded level-set distance]\label{ass:D}
There is a finite $D\ge0$ such that $\sup\{\dist(x,\Xs):x\in C,\ f(x)\le f(x_0)\}\le D$. Thus any certified upper bound may be used.
\end{assumption}

This is weaker than bounded level sets: $\Xs$ and $C$ may be unbounded. The growth condition used later implies it: Assumption~\ref{as:grow} on $U_0\ge\Delta_0$ permits the choice $D=(\Delta_0/c)^{1/p}$. Throughout this section write $D_k:=\dist(x_k,\Xs)$ and let $z_k\in\Xs$ be a nearest minimizer of $x_k$, which exists because $\Xs$ is closed, so that $D_k=\norm{x_k-z_k}$.

\begin{lemma}[small radii keep the initial level set invariant]\label{lem:inv}
Let Assumptions~\ref{as:geom}, \ref{as:conv}, \ref{as:quad}, \ref{as:att}, and~\ref{ass:D} hold with $\Delta_0>0$, and run the direct update \eqref{eq:update} with radii satisfying
\begin{equation}\label{eq:rinv}
0<r_k\le\min\Bigl\{D,\ \frac{2\Delta_0}{LD}\Bigr\}\qquad(k\ge0).
\end{equation}
Then $f(x_k)\le f(x_0)$, hence $D_k\le D$, for every $k\ge0$, and
\begin{equation}\label{eq:invrec}
\Delta_{k+1}\le\Bigl(1-\frac{r_k}D\Bigr)\Delta_k+\frac L2r_k^2 .
\end{equation}
For a constant radius $r_k\equiv r$, consequently,
\begin{equation}\label{eq:invconst}
\Delta_K\le\Bigl(1-\frac rD\Bigr)^K\Delta_0+\frac{LrD}2\qquad(K\ge0).
\end{equation}
\end{lemma}
\begin{proof}
We show by induction that $\Delta_k\le\Delta_0$. If this holds at $k$, then $f(x_k)\le f(x_0)$, so $D_k\le D$. Put $\theta:=\min\{1,r_k/D_k\}$ ($\theta=1$ if $D_k=0$). Then $\theta\norm{z_k-x_k}\le r_k$, and $r_k\le D$ gives $\theta\ge r_k/D$. The contraction \eqref{eq:contract} with $z=z_k$ yields
\begin{eqnarray*}
\Delta_{k+1} &\overset{\eqref{eq:contract}}\le& (1-\theta)\Delta_k+\frac L2r_k^2\ \le\ \Bigl(1-\frac{r_k}D\Bigr)\Delta_k+\frac L2r_k^2\\
&\le& \Delta_0-\frac{r_k}D\Delta_0+\frac L2r_k^2\ \overset{\eqref{eq:rinv}}\le\ \Delta_0,
\end{eqnarray*}
where the second inequality uses $\theta\ge r_k/D$ and $\Delta_k\ge0$, the third uses $\Delta_k\le\Delta_0$ and $1-r_k/D\ge0$, and the last is $r_k\le2\Delta_0/(LD)$. This closes the induction and proves \eqref{eq:invrec}. For a constant radius, unrolling \eqref{eq:invrec} gives
\[
\Delta_K\overset{\eqref{eq:invrec}}\le\Bigl(1-\frac rD\Bigr)^K\Delta_0+\frac L2r^2\sum_{j<K}\Bigl(1-\frac rD\Bigr)^j\le\Bigl(1-\frac rD\Bigr)^K\Delta_0+\frac L2r^2\cdot\frac Dr ,
\]
which is \eqref{eq:invconst}.
\end{proof}

The condition $r_k\le2\Delta_0/(LD)$ says that the model error $Lr_k^2/2$ of one step never exceeds the decrease $(r_k/D)\Delta_0$ that the comparator guarantees at the boundary of the level set. The gaps need not decrease, but they cannot climb above $\Delta_0$; this is all that a constant radius needs.

\begin{corollary}[direct update, constant radius: $\tilde O(1/\varepsilon)$ without transient]\label{cor:constr}
Let Assumptions~\ref{as:geom}, \ref{as:conv}, \ref{as:quad}, \ref{as:att}, and~\ref{ass:D} hold and let $0<\varepsilon\le\min\{\Delta_0,\ LD^2\}$. Run the direct update \eqref{eq:update} (full steps, possibly nonmonotone) with the constant radius and iteration count
\[
r:=\frac{\varepsilon}{LD},\qquad K:=\Bigl\lceil\frac{LD^2}{\varepsilon}\ln\frac{2\Delta_0}{\varepsilon}\Bigr\rceil
\]
(so $r\approx D\ln(2\Delta_0/\varepsilon)/K$, of order $1/K$). Then $f(x_k)\le f(x_0)$ for all $k$, and
\[
\Delta_K\ \le\ \Bigl(1-\frac{\varepsilon}{LD^2}\Bigr)^K\Delta_0+\frac\varepsilon2\ \le\ \varepsilon .
\]
The number of oracle calls to reach accuracy $\varepsilon$ is $O\bigl((LD^2/\varepsilon)\ln(2\Delta_0/\varepsilon)\bigr)$. If $\varepsilon\ge LD^2$, one step with $r=D$ already gives $\Delta_1\le LD^2/2\le\varepsilon/2$.
\end{corollary}

\begin{proof}
The target range implies $r\le D$ and $r\le2\Delta_0/(LD)$, so Lemma~\ref{lem:inv} gives
\[
\Delta_K\overset{\eqref{eq:invconst}}\le\Bigl(1-\frac{\varepsilon}{LD^2}\Bigr)^K\Delta_0+\frac\varepsilon2\le e^{-K\varepsilon/(LD^2)}\Delta_0+\frac\varepsilon2\le\varepsilon.
\]
For $\varepsilon\ge LD^2$, $r=D$ includes a nearest minimizer $z_0$, and $\Delta_1\overset{\eqref{eq:comp}}\le LD^2/2$.
\end{proof}

\begin{corollary}[calibrated radii: an all-iteration last-iterate rate]\label{cor:calibrated}
Under the hypotheses of Lemma~\ref{lem:inv}, prescribe
\begin{equation}\label{eq:calibrated}
\tau:=\max\{2,2LD^2/\Delta_0\},\qquad A:=\tau\Delta_0,\qquad r_k:=\frac{2D}{k+\tau}.
\end{equation}
Then every direct-oracle sequence stays in the initial level set and satisfies
\begin{equation}\label{eq:calibratedrate}
\Delta_k\le\frac{A}{k+\tau}=\frac{2\max\{\Delta_0,LD^2\}}{k+\tau}\qquad(k\ge0).
\end{equation}
Alternatively, under Assumptions~\ref{as:geom}, \ref{as:conv}, \ref{as:quad}, and~\ref{as:att}, if a known $D>0$ bounds $\dist(x,\Xs)$ for every $x\in C$, then the radii $r_k:=2D/(k+2)$ give
\begin{equation}\label{eq:globaldistance}
\Delta_k\le\frac{2\max\{\Delta_0,LD^2\}}{k+2}\qquad(k\ge0).
\end{equation}
This second radius rule uses only $D$, for example a known diameter bound on $C$; $L$ and $\Delta_0$ enter its guarantee, not the algorithm.
\end{corollary}
\begin{proof}
By \eqref{eq:calibrated}, $r_k\le\min\{D,\Delta_0/(LD)\}$, so Lemma~\ref{lem:inv} applies. With $n:=k+\tau\ge2$, its recurrence is
\begin{equation}\label{eq:calibratedrec}
\Delta_{k+1}\overset{\eqref{eq:invrec},\eqref{eq:calibrated}}\le
\Bigl(1-\frac2n\Bigr)\Delta_k+\frac{2LD^2}{n^2}.
\end{equation}
The induction starts at $\Delta_0=A/\tau$. If $\Delta_k\le A/n$, the nonnegative coefficient and $A\ge2LD^2$ give
\[
\Delta_{k+1}\overset{\eqref{eq:calibratedrec}}\le\Bigl(1-\frac2n\Bigr)\frac A n+\frac{2LD^2}{n^2}
\overset{\eqref{eq:calibrated}}\le\frac A n-\frac A{n^2}\le\frac A{n+1}.
\]
For the global-distance case, the comparator proof of \eqref{eq:invrec} applies at every iterate without a level-set argument. Use the same induction with $n=k+2$ and $A=2\max\{\Delta_0,LD^2\}$.
\end{proof}
Corollary~\ref{cor:calibrated} is a direct consequence of the standard scalar recurrence. It needs no monotone descent. The initial-level version uses $L,D,\Delta_0$ to calibrate the shift; replacing $\Delta_0$ by an upper bound in its denominator need not preserve invariance. The global-distance version avoids this issue. Neither version replaces the problem-constant-free guarantee below.

\paragraph{Relation to the case of a bounded feasible set} For the same update with a constant radius $\eta$ and a regularizer whose domain has diameter $D$, Kovalev~\cite[Theorem~7 and Corollary~7]{kovalev2025understanding} proves a bound of the same form as \eqref{eq:invconst}, with the error term $\frac32LD\eta$, and hence the complexity $\tilde O(LD^2/\varepsilon)$, under star-convexity and Lipschitz continuity of the gradient in the dual norm; his argument combines the one-step inequality \eqref{eq:master} with $F(x_{k+1})-F(x^\star)\le D\,\dnorm{\nabla f(x_{k+1})+\hat\nabla R_{k+1}}$. Corollary~\ref{cor:constr} is obtained by the comparator route and differs in three respects: $C$ may be unbounded, since $D$ is the level-set distance bound of Assumption~\ref{ass:D} and Lemma~\ref{lem:inv} keeps the iterates, which are not monotone, in the initial level set; only Assumption~\ref{as:quad} is used, not Lipschitz continuity of the gradient in the dual norm; and the error term is $\frac12LDr$. In the other direction, his result needs only star-convexity with respect to a minimizer, which is weaker than convexity.

Corollaries~\ref{cor:constr} and~\ref{cor:calibrated} use $D$ to set the radii. The next result removes this requirement: with radii that use no problem constant, the direct update attains $O(\log K/K)$ for all sufficiently large $K$, the constant and the onset depending on $D$.

\begin{theorem}[direct update, radii that do not use $D$]\label{thm:dfree}
Let Assumptions~\ref{as:geom}, \ref{as:conv}, \ref{as:quad}, \ref{as:att}, and~\ref{ass:D} hold with $\Delta_0>0$; the value of $D$ is not used by the algorithm. Let $a>0$ and run the direct update \eqref{eq:update} with
\[
r_k:=\frac{a\,\ell_k}{k+1},\qquad \ell_k:=\ln(k+2).
\]
Put $\bar\Delta:=\Delta_0+\frac L2\sum_{j\ge0}r_j^2<\infty$, $\bar D:=D\bar\Delta/\Delta_0$, and let $\bar k$ be the smallest index with $\ell_{\bar k}\ge2\bar D/a$ and $r_{\bar k}\le\bar D$. Then for every $K\ge\bar k$,
\[
\Delta_K\ \le\ B\,\frac{\ln(K+2)}{K+1},\qquad B:=\max\Bigl\{La\bar D,\ \frac{(\bar k+1)\bar\Delta}{\ln(\bar k+2)}\Bigr\}.
\]
\end{theorem}
\begin{proof}
\emph{Step 1: the gaps stay bounded.} Since $d_k\ge0$, Lemma~\ref{lem:comp}(b) gives $f(x_{k+1})\le f(x_k)+\frac L2r_k^2$, hence
\[
\Delta_k\ \le\ \Delta_0+\frac L2\sum_{j<k}r_j^2\ \le\ \bar\Delta\qquad(k\ge0);
\]
the series converges because $r_j^2=a^2\ln^2(j+2)/(j+1)^2$.

\emph{Step 2: the distances stay bounded.} We show $\dist(x,\Xs)\le\bar D$ for every $x\in C$ with $\Delta(x)\le\bar\Delta$. If $\Delta(x)\le\Delta_0$, this is Assumption~\ref{ass:D}. Otherwise let $z$ be a nearest minimizer of $x$, $D_x:=\norm{x-z}$, $t:=\Delta_0/\Delta(x)\in(0,1)$ and $y:=z+t(x-z)\in C$. Convexity gives $\Delta(y)\le t\Delta(x)=\Delta_0$, so $\dist(y,\Xs)\le D$. For every $z'\in\Xs$,
\[
D_x\ \le\ \norm{x-z'}\ \le\ \norm{x-y}+\norm{y-z'}=(1-t)D_x+\norm{y-z'},
\]
so $\norm{y-z'}\ge tD_x$. Taking the infimum over $z'$ gives $tD_x\le D$, i.e.\ $D_x\le D\Delta(x)/\Delta_0\le\bar D$. By Step~1, $D_k\le\bar D$ for every $k$.

\emph{Step 3: a contraction from step $\bar k$ on.} Both conditions defining $\bar k$ persist for $k\ge\bar k$, because $\ell_k$ increases and $r_k$ decreases. For such $k$ put $\theta:=\min\{1,r_k/D_k\}$ ($\theta=1$ if $D_k=0$). Then $\theta D_k\le r_k$, and $D_k,r_k\le\bar D$ give $\theta\ge r_k/\bar D$. The contraction \eqref{eq:contract} with $z=z_k$ and $\Delta_k\ge0$ yield, with $c:=a/\bar D$,
\[
\Delta_{k+1}\ \overset{\eqref{eq:contract}}\le\ (1-\theta)\Delta_k+\frac L2r_k^2\ \le\ \Bigl(1-\frac{c\,\ell_k}{k+1}\Bigr)\Delta_k+\frac{La^2\ell_k^2}{2(k+1)^2}.
\]

\emph{Step 4: induction.} We show $\Delta_k\le B\ell_k/(k+1)$ for $k\ge\bar k$. At $k=\bar k$ this follows from $\Delta_{\bar k}\le\bar\Delta$ and the definition of $B$. If it holds at $k\ge\bar k$, then $1-c\ell_k/(k+1)\ge0$ and Step~3 give
\[
\Delta_{k+1}\ \le\ \frac{B\ell_k}{k+1}-\frac{cB\ell_k^2}{(k+1)^2}+\frac{La^2\ell_k^2}{2(k+1)^2}.
\]
On the other hand, since $\ell_{k+1}\ge\ell_k$,
\begin{eqnarray*}
\frac{B\ell_{k+1}}{k+2} &\ge& \frac{B\ell_k}{k+2}\ =\ \frac{B\ell_k}{k+1}-\frac{B\ell_k}{(k+1)(k+2)}\\
&\ge& \frac{B\ell_k}{k+1}-\frac{B\ell_k}{(k+1)^2}.
\end{eqnarray*}
Comparing the two displays, $\Delta_{k+1}\le B\ell_{k+1}/(k+2)$ holds as soon as $B(c\ell_k-1)\ge\frac12La^2\ell_k$. Since $c\ell_k\ge2$ for $k\ge\bar k$, we have $c\ell_k-1\ge c\ell_k/2$, and this follows from $B\ge La^2/c=La\bar D$.
\end{proof}

The radii use neither $D$ nor $\Delta_0$, but the bound is asymptotic: the onset can be exponential in $\bar D/a$. The theorem also holds with $\bar D$ replaced by the actual distance supremum over $\{f-f_\star\le\bar\Delta\}$; the estimate $D\bar\Delta/\Delta_0$ can be loose near optimal initialization. The logarithm replaces unknown distance information: with $r_k=a/(k+1)$ the induction instead needs $a\ge2\bar D$. No lower bound shows that the logarithm is necessary. If $\Xs$ is nonempty and bounded, all sublevel sets of $f+\iota_C$ are bounded \cite[Corollary~8.7.1]{rockafellar1997convex}, so Assumption~\ref{ass:D} holds automatically.

Corollary~\ref{cor:calibrated} gives an all-iteration $O(1/K)$ rate with distance information, whereas Theorem~\ref{thm:dfree} gives an asymptotic $O(\log K/K)$ rate without using it. Whether that logarithm is necessary for radii using no problem constants remains open; this is distinct from the question without any distance-bound assumption.

\subsection{Strong convexity improves the prescribed contraction}\label{sec:strong}

Strong convexity is the first regularity condition beyond convexity, and it strengthens the comparator inequality itself: the linear comparison with a minimizer acquires the term $-\frac\mu2\norm{x-\xs}^2$, which both locates the minimizer and permits a larger radius while maintaining a gap envelope. Section~\ref{sec:growth} will use only the location, through a growth condition weaker than strong convexity, and obtain a correspondingly weaker contraction.

\begin{assumption}[strong convexity]\label{as:strong}
For some $0<\mu\le L$,
\begin{equation}\label{eq:strong}
f(y)\ge f(x)+\ip{\nabla f(x)}{y-x}+\frac\mu2\norm{y-x}^2\qquad(x,y\in C).
\end{equation}
\end{assumption}

Under Assumptions~\ref{as:geom} and~\ref{as:strong} a minimizer $\xs$ exists and is unique, so attainment need not be assumed: \eqref{eq:strong} with $x=x_0$ shows that the sublevel sets of $f$ on the closed set $C$ are bounded, and \eqref{eq:strong} with $x=\xs$, together with $\ip{\nabla f(\xs)}{y-\xs}\ge0$ for $y\in C$, gives $f(y)\ge f(\xs)+\frac\mu2\norm{y-\xs}^2$.

\begin{theorem}[stronger one-call linear convergence]\label{thm:strong}
Under Assumptions~\ref{as:geom}, \ref{as:quad}, and~\ref{as:strong}, let $U_0>0$ satisfy $\Delta_0\le U_0$. Let $\kappa:=L/\mu\ge1$ be the condition number, and prescribe
\begin{equation}\label{eq:scschedule}
(r(U),\rho):=\begin{cases}\bigl(\sqrt{2\mu U}/L,\ 1-\tfrac1\kappa\bigr), & \kappa\ge2,\\[2pt]
\bigl(\sqrt{U/(2\mu)},\ \tfrac\kappa4\bigr), & 1\le\kappa<2,\end{cases}
\qquad U_k:=\rho^kU_0,\quad r_k:=r(U_k).
\end{equation}
Then $\Delta_k\le U_k$ for every $k\ge0$ and every oracle selection, with the explicit contraction guarantee
\begin{equation}\label{eq:sccontr}
\Delta_{k+1}\le\rho U_k=U_{k+1}=\rho^{k+1}U_0\qquad(k\ge0).
\end{equation}
In particular, the contraction factor is $1-\mu/L$ when $\mu\le L/2$, and $1/4$ when $\mu=L$.
\end{theorem}
\begin{proof}
We prove, for a center $x$ with $\Delta(x)\le U$, a radius $r$ and an output $q$, the one-step envelope
\begin{equation}\label{eq:onestep}
f(q)-\fs\le\max\Bigl\{\frac L2r^2,\ U-r\sqrt{2\mu U}+\frac L2r^2\Bigr\}.
\end{equation}
Write $D_x:=\norm{x-\xs}$. If $D_x\le r$, the segment comparator \eqref{eq:segment} with $z=\xs$, $\theta=1$, and \eqref{eq:strong} give
\[
f(q)-\fs\ \overset{\eqref{eq:segment}}\le\ \Delta(x)+\ip{\nabla f(x)}{\xs-x}+\frac L2r^2\ \overset{\eqref{eq:strong}}\le\ \frac L2r^2-\frac\mu2D_x^2\ \le\ \frac L2r^2 .
\]
If $D_x>r$, use \eqref{eq:segment} with $\theta=r/D_x$, then \eqref{eq:strong}, then $\frac{\Delta}{D}+\frac\mu2D\ge\sqrt{2\mu\Delta}$:
\begin{eqnarray*}
f(q)-\fs &\overset{\eqref{eq:segment}}\le& \Delta(x)+\frac r{D_x}\ip{\nabla f(x)}{\xs-x}+\frac L2r^2\\
&\overset{\eqref{eq:strong}}\le& \Delta(x)-r\Bigl(\frac{\Delta(x)}{D_x}+\frac\mu2D_x\Bigr)+\frac L2r^2\\
&\le& \Delta(x)-r\sqrt{2\mu\Delta(x)}+\frac L2r^2 .
\end{eqnarray*}
In this case $\Delta(x)\ge\mu D_x^2/2>\mu r^2/2$ by constrained optimality of $\xs$ and \eqref{eq:strong}. The function $t\mapsto t-r\sqrt{2\mu t}+Lr^2/2$ is convex, so its maximum on $[\mu r^2/2,U]$ is attained at an endpoint: at $t=\mu r^2/2$ it equals $(L-\mu)r^2/2\le Lr^2/2$, and at $t=U$ it is the second term of \eqref{eq:onestep}. This proves \eqref{eq:onestep}.

Now insert \eqref{eq:scschedule}. For $\kappa\ge2$, $r=\sqrt{2\mu U}/L$ and
\[
\frac L2r^2\overset{\eqref{eq:scschedule}}=\frac U\kappa,\qquad U-r\sqrt{2\mu U}+\frac L2r^2\overset{\eqref{eq:scschedule}}=U-\frac{2U}\kappa+\frac U\kappa=\Bigl(1-\frac1\kappa\Bigr)U .
\]
For $1\le\kappa<2$, $r=\sqrt{U/(2\mu)}$ and both terms equal $\kappa U/4$. In either case the bound is $\rho U$, and applying it with $U=U_k$ proves \eqref{eq:sccontr} by induction.
\end{proof}

Equation~\eqref{eq:sccontr} contracts the prescribed envelope $U_k$, not necessarily the actual gap; the same holds for every prescribed schedule below. The argument allows a constrained minimizer with $\nabla f(\xs)\ne0$. The factor, $1-\mu/L$ for $\mu\le L/2$ and $L/(4\mu)$ for $L/2<\mu\le L$, improves the $1-\mu/(4L)$ of the quadratic-growth specialization of Theorem~\ref{thm:growth} by using the full inequality in Assumption~\ref{as:strong}. An initial envelope is available without knowledge of $\fs$: $U_0:=\dnorm{\nabla f(x_0)}^2/(2\mu)$ satisfies $\Delta_0\le U_0$, Indeed, with $D=\norm{x_0-\xs}$,

\[\Delta_0\overset{\eqref{eq:strong}}\le\ip{\nabla f(x_0)}{x_0-\xs}-\frac\mu2D^2\le\dnorm{\nabla f(x_0)}D-\frac\mu2D^2\le\frac{\dnorm{\nabla f(x_0)}^2}{2\mu}.\]

If this upper bound is zero, stop: $x_0$ is optimal. Otherwise it is a valid positive $U_0$.

\subsection{Growth calibrates radii and controls total movement}\label{sec:growth}

Growth converts a gap bound $U$ into a distance bound $R=(U/c)^{1/p}$. Taking radius $r=\theta R$ permits a fraction $\theta$ of the comparison toward a minimizer. We choose $\theta$ so that the model error $Lr^2/2$ consumes at most half of the guaranteed decrease $\theta U$.

\begin{assumption}[power growth]\label{as:grow}
For fixed $c>0$, $p\ge1$, and $U_0>0$ with $\Delta_0\le U_0$,
\begin{equation}\label{eq:growth}
f(x)-\fs\ge c\dist(x,\Xs)^p\qquad(x\in C,\ f(x)-\fs\le U_0).
\end{equation}
\end{assumption}

\begin{theorem}[prescribed growth-based radii]\label{thm:growth}
Under Assumptions~\ref{as:geom}, \ref{as:conv}, \ref{as:quad}, \ref{as:att}, and~\ref{as:grow}, define
\begin{equation}\label{eq:growthrec}
\theta_k:=\min\Bigl\{1,\frac{c^{2/p}}LU_k^{1-2/p}\Bigr\},\qquad r_k:=\theta_k(U_k/c)^{1/p},\qquad U_{k+1}:=(1-\theta_k)U_k+\frac L2r_k^2 .
\end{equation}
Every direct-oracle sequence satisfies $\Delta_k\le U_k\to0$. Moreover, it converges to a minimizer $\xs\in\Xs$ selected by the sequence, with
\begin{equation}\label{eq:tail}
\norm{x_k-\xs}\le\sum_{j\ge k}\norm{x_{j+1}-x_j}\le\sum_{j\ge k}r_j\le2p(U_k/c)^{1/p}.
\end{equation}
\end{theorem}
\begin{proof}
Suppose $\Delta_k\le U_k\le U_0$. A nearest minimizer $z_k$ exists because $\Xs$ is closed, and Assumption~\ref{as:grow} gives $\norm{x_k-z_k}\le(U_k/c)^{1/p}$, so $\theta_k\norm{z_k-x_k}\le r_k$. The contraction \eqref{eq:contract} with $z=z_k$, together with $1-\theta_k\ge0$ and $\Delta_k\le U_k$, gives
\[
\Delta_{k+1}\ \overset{\eqref{eq:contract}}\le\ (1-\theta_k)\Delta_k+\frac L2r_k^2\ \le\ (1-\theta_k)U_k+\frac L2r_k^2\overset{\eqref{eq:growthrec}}=U_{k+1}.
\]
Moreover, $\theta_k\le c^{2/p}U_k^{1-2/p}/L$ means $L\theta_k(U_k/c)^{2/p}\le U_k$, so $\frac L2r_k^2=\frac L2\theta_k^2(U_k/c)^{2/p}\le\frac{\theta_k}2U_k$ and
\begin{equation}\label{eq:halfdec}
0<U_{k+1}\le(1-\theta_k/2)U_k<U_k .
\end{equation}
Thus the induction stays in the region of Assumption~\ref{as:grow}. If $U_k$ had a positive limit, so would $\theta_k$, contradicting \eqref{eq:halfdec}; hence $U_k\to0$.

For the tail bound, concavity of $t\mapsto t^{1/p}$ and \eqref{eq:halfdec} give
\[
U_k^{1/p}-U_{k+1}^{1/p}\ \ge\ \frac1pU_k^{1/p-1}(U_k-U_{k+1})\ \overset{\eqref{eq:halfdec}}\ge\ \frac{\theta_k}{2p}U_k^{1/p}\overset{\eqref{eq:growthrec}}=\frac{c^{1/p}}{2p}\,r_k .
\]
Summing over $j\ge k$ and using $U_j\to0$ gives $\sum_{j\ge k}r_j\le2p(U_k/c)^{1/p}$, and $\norm{x_{j+1}-x_j}\le r_j$ gives the remaining inequalities in \eqref{eq:tail}. Hence $(x_k)$ is Cauchy; its limit lies in $C$ because $C$ is closed, and in $\Xs$ because $\Delta_k\to0$ and $f$ is continuous.
\end{proof}

\paragraph{The three growth regimes} Under the hypotheses of Theorem~\ref{thm:growth}, the scalar envelope has explicit rates. If $p=2$, then $\theta:=\min\{1,c/L\}$ is constant and $U_k\le(1-\theta/2)^kU_0$. Quadratic growth with $c=\mu/2$ and $\mu\le L$ gives the factor $1-\mu/(4L)$; Theorem~\ref{thm:strong}, which uses the full strong-convexity inequality and not only the location of the minimizer, improves it to $1-\mu/L$ for $\mu\le L/2$ and to $L/(4\mu)$ for $L/2<\mu\le L$. Equation~\eqref{eq:tail} additionally gives convergence to a selected minimizer, even when it is not unique.

If $p>2$, eventually $\theta_k<1$. Put $\xi:=1-2/p>0$ and $\kappa_p:=c^{2/p}/(2L)$. From any such index $m$ onward,
\begin{equation}\label{eq:pbig}
U_{k+1}=U_k-\kappa_pU_k^{1+\xi},\qquad U_k\le\bigl(U_m^{-\xi}+\xi\kappa_p(k-m)\bigr)^{-1/\xi}=O(k^{-p/(p-2)}).
\end{equation}
Indeed, once $\theta_k<1$ the envelope reads $U_{k+1}=U_k-\kappa_pU_k^{1+\xi}$ with $v:=\kappa_pU_k^{\xi}=\theta_k/2<\tfrac12$, and since $(1-v)^{-\xi}\ge1+\xi v$ for $0\le v<1$, one gets $U_{k+1}^{-\xi}\ge U_k^{-\xi}+\xi\kappa_p$; summing from $m$ gives \eqref{eq:pbig}. The endpoint error in \eqref{eq:tail} is then $O(k^{-1/(p-2)})$.

If $1\le p<2$, eventually $\theta_k=1$, and
\begin{equation}\label{eq:psmall}
U_{k+1}=\frac L{2c^{2/p}}U_k^{2/p}.
\end{equation}
This is a superlinear envelope, quadratic when $p=1$.

\section{Nonconvex objectives}\label{sec:nonconvex}

Without convexity, the comparison of Section~\ref{sec:onestep} no longer turns into a bound on objective values, and the analysis relies instead on the certificate at the output (Lemma~\ref{lem:cert}), which needs the gradient Lipschitz condition in Assumption~\ref{as:dual}. Summing the descent inequality and selecting the output with the smallest observed model decrease gives $O(K^{-1/2})$ stationarity, measured by the residual $s$ in the dual norm (Section~\ref{sec:statsub}). When $s$ controls the objective gap through a constrained Polyak--{\L}ojasiewicz inequality, prescribed geometric radii give linear convergence (Section~\ref{sec:plsub}). Finally, an observable radius-halving test removes the PL modulus and the initial gap bound from the algorithm, at the cost of extra oracle calls, which are counted exactly when the search is warm-started (Section~\ref{sec:adaptive}).

\subsection{Stationarity of a selected output}\label{sec:statsub}

\begin{theorem}[observable stationarity bound]\label{thm:stat}
Under Assumptions~\ref{as:geom} and~\ref{as:dual}, run $K\ge1$ updates with radius $r$. Return $\hat x_K:=x_{i+1}$ for any $i\in\argmin_{0\le k<K}d_k$, where $d_k:=-\ip{\nabla f(x_k)}{x_{k+1}-x_k}$. Then
\begin{equation}\label{eq:statrate}
s(\hat x_K)\le\frac{f(x_0)-\fs}{Kr}+\frac{3L}2r .
\end{equation}
In particular, $r=a/\sqrt K$ gives an $O(K^{-1/2})$ bound for any $a>0$.
\end{theorem}
\begin{proof}
Summing Lemma~\ref{lem:comp}(b) over $k<K$ gives
\[
\sum_{k<K}d_k\ \le\ f(x_0)-f(x_K)+\frac L2Kr^2\ \le\ f(x_0)-\fs+\frac L2Kr^2 .
\]
For $i\in\argmin_kd_k$, $d_i$ is at most the average of the $d_k$, so the certificate \eqref{eq:cert} at $x_{i+1}$ gives
\[
s(\hat x_K)\ \overset{\eqref{eq:cert}}\le\ Lr+\frac{d_i}r\ \le\ Lr+\frac1{Kr}\sum_{k<K}d_k\ \le\ \frac{f(x_0)-\fs}{Kr}+\frac{3L}2r .
\]
\end{proof}

The bound \eqref{eq:statrate} is that of Kovalev~\cite[Theorem~1]{kovalev2025understanding}, who proves it for $\min_{1\le k\le K}\dnorm{\nabla f(x_k)+\hat\nabla R_k}$ and a general regularizer $R$; for $R=\iota_C$ this quantity is at least $\min_{1\le k\le K}s(x_k)$. Theorem~\ref{thm:stat} adds the rule that selects the output: it needs only the observed linear decreases, not the stationarity measure itself. If the initial gap is known and positive, $r=\sqrt{2(f(x_0)-\fs)/(3LK)}$ gives $s(\hat x_K)\le\sqrt{6L(f(x_0)-\fs)/K}$.

\subsection{A constrained PL inequality gives linear convergence}\label{sec:plsub}

The certificate becomes an objective guarantee if a small residual forces a small gap. With constraints, the gradient itself may be large while its feasible descent is small, so we impose this condition on $s(x)$ from \eqref{eq:normalcone}.

\begin{assumption}[Polyak--{\L}ojasiewicz inequality for the normal-cone residual]\label{as:pl}
We call the following condition constrained PL. For some $\mu>0$ and finite $\fs=\inf_Cf$,
\begin{equation}\label{eq:pl}
\frac12s(x)^2\ge\mu\bigl(f(x)-\fs\bigr)\qquad(x\in C).
\end{equation}
\end{assumption}

If $C=\E$, then $\NC(x)=\{0\}$ and this is the ordinary PL inequality in the chosen dual norm. For general $C$, it instead measures the residual after allowing a constraint normal to cancel the gradient. In the Euclidean norm, \eqref{eq:pl} is weaker than the PL condition on the projected-gradient mapping $G_\gamma(x):=(x-P_C(x-\gamma\nabla f(x)))/\gamma$, where $P_C$ is the projection onto $C$, that is used by Richt{\'a}rik et~al.~\cite{richtarik2026broximal}. Since $x=P_C(x+\gamma n)$ for every $n\in\NC(x)$, nonexpansiveness of $P_C$ gives $\norm{G_\gamma(x)}_2\le\norm{\nabla f(x)+n}_2$, hence $\norm{G_\gamma(x)}_2\le s(x)$, so the projected-gradient condition implies \eqref{eq:pl} with the same constant. The converse fails: for $f(x)=x$ on $C=[0,1]$ one has $s(x)=1$ for $x>0$ and $s(0)=0$, so \eqref{eq:pl} holds with $\mu=\frac12$, whereas $G_1(x)=x$ and $G_1(x)^2/(2\Delta(x))=x/2\to0$ as $x\downarrow0$. Generalized PL conditions already yield linear convergence of proximal methods \cite{karimi2016linear}; the result below derives a prescribed-radius consequence for BroxGD from the one-step inequality \eqref{eq:master}.

\begin{theorem}[linear convergence under constrained PL]\label{thm:pl}
Under Assumptions~\ref{as:geom}, \ref{as:dual}, and~\ref{as:pl}, choose $U_0>0$ with $\Delta_0\le U_0$ and prescribe
\begin{equation}\label{eq:plschedule}
\rho:=\frac{3L}{3L+\mu},\qquad U_k:=\rho^kU_0,\qquad r_k:=\frac{\sqrt{2\mu\rho U_k}}{3L}.
\end{equation}
Every exact BroxGD sequence satisfies
\begin{equation}\label{eq:plrate}
\Delta_{k+1}\le\rho U_k=U_{k+1},\qquad \Delta_k\le\rho^kU_0 .
\end{equation}
\end{theorem}
\begin{proof}
Consider an update from $x$ to $q$ with radius $r$. Assumption~\ref{as:pl} gives $s(q)\ge\sqrt{2\mu\Delta(q)}$, and inserting this into \eqref{eq:master} yields
\begin{equation}\label{eq:plkey}
\Delta(q)+r\sqrt{2\mu\Delta(q)}\ \overset{\eqref{eq:pl}}\le\ \Delta(q)+r\,s(q)\ \overset{\eqref{eq:master}}\le\ \Delta(x)+\frac{3L}2r^2 .
\end{equation}
Suppose $\Delta(x)\le U$ and $r=\sqrt{2\mu\rho U}/(3L)$ as in \eqref{eq:plschedule}. Since $1-\rho=\mu\rho/(3L)$,
\[
\frac{3L}2r^2\overset{\eqref{eq:plschedule}}=\frac{\mu\rho U}{3L}=(1-\rho)U,\qquad r\sqrt{2\mu\rho U}\overset{\eqref{eq:plschedule}}=\frac{2\mu\rho U}{3L}=2(1-\rho)U ,
\]
so that
\[
\rho U+r\sqrt{2\mu\rho U}=U+(1-\rho)U=U+\frac{3L}2r^2\ \ge\ \Delta(x)+\frac{3L}2r^2 .
\]
The left side of \eqref{eq:plkey} is strictly increasing in $\Delta(q)\ge0$, so comparing the last display with \eqref{eq:plkey} gives $\Delta(q)\le\rho U$. Induction with $U=U_k$ proves \eqref{eq:plrate}.
\end{proof}

For $0<\varepsilon<U_0$, \eqref{eq:plrate} gives $\Delta_K\le\varepsilon$ whenever
\[
K\ge\Bigl\lceil\frac{\log(U_0/\varepsilon)}{\log(1+\mu/(3L))}\Bigr\rceil .
\]
As under strong convexity, $U_0:=\dnorm{\nabla f(x_0)}^2/(2\mu)$ is a valid initial envelope, because $\Delta_0\le s(x_0)^2/(2\mu)$ by \eqref{eq:pl} and $s(x_0)\le\dnorm{\nabla f(x_0)}$; no normal-cone computation is needed. If this choice gives $U_0=0$, stop at the already optimal $x_0$ rather than generate zero radii.

\subsection{Adaptive radii without a supplied PL constant}\label{sec:adaptive}

Prescribed envelopes need quantitative regularity information. We can remove the PL constant and the initial gap bound from the algorithm by testing the observable linear decrease, at the cost of extra oracle calls. The construction assumes a known valid gradient Lipschitz constant in the dual norm $L$. 

At each center $x$, start at a chosen cap $\hat r>0$ and query the exact BroxGD oracle at radii $r=\hat r,\hat r/2,\hat r/4,\dots$. Write $d(r)=-\ip{\nabla f(x)}{q(r)-x}$. If $d(r)=0$, stop and return $x$. Otherwise accept the first radius satisfying
\begin{equation}\label{eq:halving}
d(r)\ge Lr^2,
\end{equation}
and set the next iterate to $q(r)$. Restart from $\hat r$ at the next center. The gradient is evaluated once per center and reused during its radius search. The test needs only $L$, $\hat r$, and oracle output values, not $\mu$, $\fs$, or a gap bound.

\begin{theorem}[radius halving and constrained PL]\label{thm:halving}
Under Assumptions~\ref{as:geom} and~\ref{as:dual}, the search terminates after finitely many calls at every nonstationary center. Every accepted update decreases the objective. Put $A_{\hat r}=L\hat r^2/2$. If Assumption~\ref{as:pl} also holds, an early stop is globally optimal, and every accepted update satisfies
\begin{equation}\label{eq:halvingdec}
\Delta_k-\Delta_{k+1}\ge\min\Bigl\{A_{\hat r},\ \frac{\mu\Delta_{k+1}}{25L}\Bigr\}.
\end{equation}
Updates accepted at $r=\hat r$ decrease the gap by at least $A_{\hat r}$; all other accepted updates satisfy $\Delta_{k+1}\le\hat\rho\Delta_k$, where $\hat\rho=25L/(25L+\mu)$. For $0<\varepsilon<\Delta_0$, reaching $\Delta_k\le\varepsilon$ takes at most
\begin{equation}\label{eq:halvingbudget}
\Bigl\lfloor\frac{\Delta_0}{A_{\hat r}}\Bigr\rfloor+\Bigl\lceil\frac{\log(\Delta_0/\varepsilon)}{\log(1+\mu/(25L))}\Bigr\rceil
\end{equation}
accepted updates, or an earlier optimal stop. This is not a bound on total oracle calls.
\end{theorem}
\begin{proof}
\emph{Finite search and early stop.} At a nonstationary center, $-\nabla f(x)\notin\NC(x)$, so there is $y\in C$ with $\gamma:=-\ip{\nabla f(x)}{y-x}>0$. For $0<r\le\tau:=\norm{y-x}$, comparison with $x+(r/\tau)(y-x)$ gives $d(r)\ge\gamma r/\tau$, so \eqref{eq:halving} holds as soon as $r\le\gamma/(L\tau)$, and the search is finite. If $d(r)=0$ for some tried radius $r>0$, then $s(x)=0$ by Lemma~\ref{lem:dr}(iii); by the same lemma this can only happen at the first radius $\hat r$ tried at $x$, so an early stop costs exactly one oracle call and occurs only at a stationary center. Under Assumption~\ref{as:pl}, $s(x)=0$ makes $x$ globally optimal.

\emph{Decrease of an accepted update.} By Lemma~\ref{lem:comp}(b) and the acceptance test $d\ge Lr^2$,
\[
\Delta_k-\Delta_{k+1}\ \ge\ d-\frac L2r^2\ \overset{\eqref{eq:halving}}\ge\ \frac d2\ \overset{\eqref{eq:halving}}\ge\ \frac L2r^2 .
\]
If $r=\hat r$, this is at least $A_{\hat r}$. Otherwise the radius $2r$ was rejected, and since $d$ is nondecreasing, $d(r)\le d(2r)<4Lr^2$. Then \eqref{eq:cert} gives $s(x_{k+1})<5Lr$, and
\[
\Delta_k-\Delta_{k+1}\ \ge\ \frac L2r^2\ >\ \frac{s(x_{k+1})^2}{50L}\ \overset{\eqref{eq:pl}}\ge\ \frac{\mu\Delta_{k+1}}{25L} ,
\]
which proves \eqref{eq:halvingdec} and $\Delta_{k+1}\le\hat\rho\Delta_k$ for uncapped steps.

\emph{Budget.} Capped steps decrease the gap by at least $A_{\hat r}$ and the gap is nonnegative, so there are at most $\lfloor\Delta_0/A_{\hat r}\rfloor$ of them. All steps decrease the gap and each uncapped step multiplies it by at most $\hat\rho=1/(1+\mu/(25L))$, so after $m$ uncapped steps the gap is at most $\hat\rho^m\Delta_0$, which is at most $\varepsilon$ once $m\ge\log(\Delta_0/\varepsilon)/\log(1+\mu/(25L))$, wherever the capped steps occur.
\end{proof}

The method eventually has only uncapped updates unless it stops, so it converges linearly in accepted updates under constrained PL. Every search is finite, but the proof supplies no uniform bound on rejected calls, so accepted-update counts are not comparable with the one-call prescribed schedules; the warm start below repairs this. Radius backtracking for local linear solves goes back to successive linear programming \cite{ferris1996linear}.

\paragraph{Warm-started halving} The rejected calls can be counted, and then bounded, if each search starts at twice the previously accepted radius instead of at the fixed cap $\hat r$. This is the amortized count of backtracking steps of the universal gradient method of Nesterov~\cite{nesterov2015universal}; the Armijo radius search of Ferris and Zavriev~\cite{ferris1996linear} is likewise adaptive.
\begin{algo}[warm-started radius halving]\label{alg:warm}
Given $\hat r_0>0$. At center $x_k$ with cap $\hat r_k$, test $r=\hat r_k,\hat r_k/2,\dots$ in turn and stop if $d(r)=0$; otherwise set
\begin{equation}\label{eq:warmradius}
\begin{aligned}
r_k&:=\max\bigl\{r\in\{\hat r_k2^{-j}:j=0,1,\dots\}:\ d(r)\ge Lr^2\bigr\},\\
x_{k+1}&:=q(r_k),\qquad \hat r_{k+1}:=2r_k .
\end{aligned}
\end{equation}
\end{algo}
\begin{theorem}[total oracle calls of warm-started halving]\label{thm:warm}
Let Assumptions~\ref{as:geom}, \ref{as:dual}, and~\ref{as:pl} hold. Then
\begin{enumerate}[label=(\alph*),leftmargin=1.6em]
\item if Algorithm~\ref{alg:warm} performs at least $K$ updates, the total number of oracle calls during the first $K$ updates is exactly $2K-1+\log_2(\hat r_0/r_{K-1})$; if instead it stops at a center $x_k$ with $k<K$ because $d(r)=0$, then $s(x_k)=0$, the stop occurs at the first radius tried at $x_k$, and the run has used $2k+\log_2(\hat r_0/r_{k-1})$ calls in total ($1$ call if $k=0$);
\item let $k_{\rm cap}$ be the number of initial consecutive ``capped'' steps ($r_k=\hat r_k$); then $k_{\rm cap}\le\log_4\bigl(1+6\Delta_0/(L\hat r_0^2)\bigr)$, and for every $k\ge k_{\rm cap}$ at which an update is made, $\Delta_{k+1}\le\hat\rho\,\Delta_k$ with $\hat\rho=25L/(25L+\mu)$. Consequently $\Delta_K\le\hat\rho^{\,K-k_{\rm cap}}\Delta_0$ whenever $K$ updates are made, and $\Delta_k=0$ if the run stops at $x_k$.
\end{enumerate}
\end{theorem}
\begin{proof}
\emph{(a)} The center $x_k$ uses $1+j_k$ calls, where $j_k=\log_2(\hat r_k/r_k)\in\{0,1,2,\dots\}$ is the number of rejections. For $k\ge1$, $\hat r_k=2r_{k-1}$, so $j_k=1+\log_2(r_{k-1}/r_k)$ and
\[
\sum_{k<K}(1+j_k)=K+\log_2\frac{\hat r_0}{r_0}+\sum_{k=1}^{K-1}\Bigl(1+\log_2\frac{r_{k-1}}{r_k}\Bigr)=2K-1+\log_2\frac{\hat r_0}{r_{K-1}} .
\]

\emph{Early stop.} Lemma~\ref{lem:dr}(iii) implies that a stop can occur only on the first query at a stationary center. Add one call to the identity for the preceding $k$ updates, or count one call if $k=0$.

\emph{(b)} As in the proof of Theorem~\ref{thm:halving}, every accepted step satisfies $\Delta_k-\Delta_{k+1}\ge Lr_k^2/2$, so the gaps decrease.

\emph{Uncapped steps} ($j_k\ge1$). The radius $2r_k$ was rejected, so $d_k\le d(2r_k)<4Lr_k^2$, and \eqref{eq:cert} and \eqref{eq:pl} give
\[
s(x_{k+1})\overset{\eqref{eq:cert}}<5Lr_k,\qquad \Delta_{k+1}\ \overset{\eqref{eq:pl}}\le\ \frac{s(x_{k+1})^2}{2\mu}\ <\ \frac{25L^2r_k^2}{2\mu} .
\]
Hence $\Delta_k-\Delta_{k+1}\ge Lr_k^2/2\ge\mu\Delta_{k+1}/(25L)$, i.e.\ $\Delta_{k+1}\le\hat\rho\Delta_k$.

\emph{Capped steps after an uncapped one} ($k>k_{\rm cap}$, $r_k=\hat r_k$). Let $m<k$ be the most recent uncapped step, so that the steps $m+1,\dots,k$ are capped and $r_k=2^{k-m}r_m\ge2r_m$. By the previous paragraph, $\Delta_{m+1}<25L^2r_m^2/(2\mu)$, and the gaps decrease, so
\[
\Delta_k-\Delta_{k+1}\ \ge\ \frac L2r_k^2\ \ge\ 2Lr_m^2\ \ge\ \frac{4\mu}{25L}\Delta_{m+1}\ \ge\ \frac{4\mu}{25L}\Delta_k .
\]
Thus $\Delta_{k+1}\le(1-\frac{4\mu}{25L})\Delta_k\le\hat\rho\Delta_k$, the last inequality being equivalent to $4(25L+\mu)\ge25L$. If $1-4\mu/(25L)<0$, the preceding decrease bound rules out this capped-step configuration at a positive gap; no restriction $\mu\le L$ is needed.

\emph{Initial capped steps.} For $k<k_{\rm cap}$, $r_k=2^k\hat r_0$ and the decrease is at least $L\hat r_0^2\,4^k/2$; since the total decrease is at most $\Delta_0$,
\[
\frac{L\hat r_0^2}2\cdot\frac{4^{k_{\rm cap}}-1}3=\frac{L\hat r_0^2}2\sum_{k<k_{\rm cap}}4^k\ \le\ \Delta_0 ,
\]
which gives the bound on $k_{\rm cap}$. If $d(r)=0$ occurs, that point is globally optimal by \eqref{eq:pl}.
\end{proof}

The warm start preserves the linear rate after the initial capped steps. Its exact call identity still contains a final-radius term, so an accuracy bound needs a lower bound on accepted radii before the target is reached. The next corollary supplies it. Write $\log_2^+t:=\max\{\log_2t,0\}$.

\begin{corollary}[oracle calls to reach a prescribed accuracy]\label{cor:warmcalls}
Let Assumptions~\ref{as:geom}, \ref{as:dual}, and~\ref{as:pl} hold, let $0<\varepsilon<\Delta_0$, and let $T\ge1$ be the first index with $\Delta_T\le\varepsilon$ in a run of Algorithm~\ref{alg:warm}. Put $C_{\rm rej}:=2L^2/\mu+3L/2$. Then
\begin{equation}\label{eq:radiusfloor}
r_k\ \ge\ \min\Bigl\{\hat r_0,\ \frac12\sqrt{\frac{\varepsilon}{C_{\rm rej}}}\Bigr\}\qquad(k<T),
\end{equation}
and the number $N_T$ of oracle calls during the first $T$ updates satisfies
\begin{equation}\label{eq:warmcalls}
N_T\ \le\ 2T-1+\log_2^+\Bigl(2\hat r_0\sqrt{\frac{C_{\rm rej}}{\varepsilon}}\Bigr),\qquad T\ \le\ k_{\rm cap}+\Bigl\lceil\frac{\log(\Delta_0/\varepsilon)}{\log(1+\mu/(25L))}\Bigr\rceil .
\end{equation}
\end{corollary}
\begin{proof}
The index $T$ exists and satisfies the second bound in \eqref{eq:warmcalls} by Theorem~\ref{thm:warm}(b): either the run stops at some $x_k$ with $\Delta_k=0$, or $\Delta_K\le\hat\rho^{\,K-k_{\rm cap}}\Delta_0$ for all $K$, and $1/\hat\rho=1+\mu/(25L)$.

\emph{Rejected radii.} Let $x$ be a center with $\Delta(x)>\varepsilon$ and let $r$ be a radius rejected at $x$, with output $q$ and $d=d(r)<Lr^2$. By \eqref{eq:cert} and \eqref{eq:pl},
\[
s(q)\ \overset{\eqref{eq:cert}}\le\ Lr+\frac dr\ <\ 2Lr,\qquad \Delta(q)\ \overset{\eqref{eq:pl}}\le\ \frac{s(q)^2}{2\mu}\ <\ \frac{2L^2}{\mu}r^2 .
\]
Integration along the segment $[x,q]\subseteq C$ and Assumption~\ref{as:dual} give $f(q)\ge f(x)+\ip{\nabla f(x)}{q-x}-\frac L2\norm{q-x}^2$, that is, $\Delta(x)\le\Delta(q)+d+\frac L2r^2$. Hence $\Delta(x)<C_{\rm rej}r^2$, and every radius rejected at $x$ exceeds $\sqrt{\varepsilon/C_{\rm rej}}$.

\emph{Accepted radii.} We prove \eqref{eq:radiusfloor} by induction on $k<T$; note that $\Delta_k>\varepsilon$ for these $k$. If $r_k<\hat r_k$, the radius $2r_k$ was rejected at $x_k$, so $2r_k>\sqrt{\varepsilon/C_{\rm rej}}$. If $r_k=\hat r_k$, then either $k=0$ and $r_0=\hat r_0$, or $k\ge1$ and $r_k=2r_{k-1}\ge r_{k-1}$, to which the induction hypothesis applies.

\emph{Calls.} By Theorem~\ref{thm:warm}(a) and \eqref{eq:radiusfloor} at $k=T-1$,

\[N_T=2T-1+\log_2\frac{\hat r_0}{r_{T-1}}\overset{\eqref{eq:radiusfloor}}\le2T-1+\log_2^+\Bigl(2\hat r_0\sqrt{\frac{C_{\rm rej}}\varepsilon}\Bigr).\]

\end{proof}

Thus the gap first falls below $\varepsilon$ after $O(\log(1/\varepsilon))$ oracle calls. This is a bound on a hitting time: the algorithm does not observe $\Delta_k$, and without $\mu$ and $\fs$ it cannot recognize that the accuracy has been reached, except through the stopping rule $d(r)=0$; what it observes is the stationarity bound $s(x_{k+1})\le Lr_k+d_k/r_k$ of \eqref{eq:cert}. The proof uses \eqref{eq:pl} at rejected outputs, which may leave the level set of the current iterate; this is why \eqref{eq:pl} is assumed on all of $C$.

\section{How the regularity assumptions are related}\label{sec:related}

The assumptions describe overlapping classes. Under Assumptions~\ref{as:geom}, \ref{as:conv}, and~\ref{as:att}, strong convexity with modulus $\mu$ implies quadratic growth with $c=\mu/2$. More generally, quadratic growth with constant $c>0$ implies constrained PL with modulus $c/2$: with $D=\dist(x,\Xs)$, convexity gives $\Delta(x)\le s(x)D$. Indeed, if $z\in\Xs$ and $n\in\NC(x)$, then $\ip n{z-x}\le0$ gives
\[
\Delta(x)\le\ip{\nabla f(x)}{x-z}\le\ip{\nabla f(x)+n}{x-z}\le\dnorm{\nabla f(x)+n}\,\norm{x-z},
\]
and taking the infimum over $n$ and $z$ yields the claim. Consequently
\begin{equation}\label{eq:growth2pl}
cD^2\overset{\eqref{eq:growth}}\le\Delta(x)\le s(x)D,\qquad \frac12s(x)^2\ge\frac c2\Delta(x).
\end{equation}
These implications hold on the same region as the growth hypothesis; a local growth statement does not establish Assumption~\ref{as:pl} globally. Full strong convexity gives the sharper PL modulus $\mu$: under Assumption~\ref{as:strong} the same comparison gives $\Delta(x)\le s(x)D-\mu D^2/2$, and maximizing over $D\ge0$ gives $\tfrac12s(x)^2\ge\mu\Delta(x)$.

For $p=1$, Assumption~\ref{as:grow} is weak sharpness of the solution set in the sense of Burke and Ferris~\cite{burke1993weak}, restricted to a sublevel set. On a sublevel set with $0<\Delta(x)\le U$, growth of order $p=1$ with constant $\alpha$ gives quadratic growth with $c=\alpha^2/U$ and a PL modulus $\alpha^2/(2U)$. It does not imply global quadratic growth on an unbounded set: $f(x)=x$ on $[0,\infty)$ is weakly sharp but has no positive global quadratic-growth constant.

For convex objectives with quadratic growth (Assumption~\ref{as:grow} with $p=2$), Theorems~\ref{thm:growth} and \ref{thm:pl} therefore overlap. The former uses only the upper model and proves finite movement; the latter additionally assumes the gradient Lipschitz condition, but also covers nonconvex PL objectives for which growth alone would not justify the comparator proof. In fact its geometric radii also give finite movement: $\sum_{j\ge k}r_j=r_k/(1-\sqrt\rho)$, so completeness and continuity imply convergence to a minimizer.

\begin{table}[t]

\caption{Information used to generate radii. Constants appearing only in the guarantee need not be algorithm inputs. A supplied $U_0$ bounds the initial gap.}\label{tab:inputs}
\centering\small
\begin{tabular}{@{}p{.43\textwidth}p{.53\textwidth}@{}}\toprule
Schedule & Inputs beyond gradients and local-oracle outputs\\\midrule
Corollaries~\ref{cor:value}, \ref{cor:fixed}; Theorem~\ref{thm:dfree} & Chosen scale/exponent or horizon; no problem constant.\\
Corollary~\ref{cor:constr} & $L,D,\varepsilon$ for the radius; a gap bound for the budget.\\
Corollary~\ref{cor:calibrated}, bounded level-set / global distance & $L,D,\Delta_0$ for the shift / global distance bound $D$ only.\\
Theorems~\ref{thm:strong}, \ref{thm:growth}, \ref{thm:pl} & $L,\mu,U_0$ or $L,c,p,U_0$, as applicable.\\
Theorem~\ref{thm:stat}, $r=a/\sqrt K$ & Scale $a$ and horizon $K$.\\
Theorems~\ref{thm:halving}, \ref{thm:warm} & Valid $L$ and initial cap; neither $\mu$ nor $U_0$.\\\bottomrule
\end{tabular}
\end{table}
Table~\ref{tab:inputs} separates radius computation from performance certification. In particular, a rule can generate iterates without knowing $f_\star$ while its objective-accuracy stopping budget still needs an error bound. An oracle call is the unit of work here; these guarantees do not assert inexpensive local solves or comparable runtimes across geometries.

\section{Conclusions}\label{sec:concl}

In a general norm the Fej\'er-type decrease of a squared distance is not available, and all results of this paper rest instead on two inequalities about one step: the comparison of the linear model at the output with any other point of the local ball, and a computable bound on the stationarity measure $s$ at the output. In our analysis the price of the missing distance decrease is the control of the distance to the solution set: without it the direct update gives $f(x_k)\to\inf_Cf$ and a last-iterate $O(K^{-1/3})$, whereas distance control gives an all-iteration $O(1/K)$ bound with calibrated radii, or an asymptotic $O(\log K/K)$ bound with radii using no problem constant. Strong convexity and power growth give faster envelopes, and warm-started halving has logarithmic total oracle complexity under global constrained PL.

One question remains open: for the direct update in a general norm and without a distance bound, we do not know whether the gap between the proved $O(K^{-1/3})$ and the Euclidean $O(1/K)$ is real.

\section*{Declarations}
OpenAI Codex assisted with manuscript revision, proof checking, and the derivation and exposition of Corollary~\ref{cor:calibrated}. The authors assume responsibility for all content.

\appendix
\section{Proof of the statement in Remark~\ref{rem:ce}}\label{app:ce}

The translation argument avoids all nondifferentiability lines, ensuring that both the observed subgradients and the radii are independent of the selection.

Fix a radius rule as in Remark~\ref{rem:ce} and $x_0$. The function $f_c=f(\cdot-c)$ fails to be differentiable exactly on the three lines $x_1=c_1$, $x_2=c_2$ and $x_1+x_2=c_1+c_2$; at every other point its gradient is one of the eight vectors $(\sigma+\varepsilon a,\sigma+\varepsilon b)$, $\sigma,a,b\in\{\pm1\}$, and the oracle output at radius $r$ is $x-r\,\sigma(1,1)$. Let $P$ be the set of all points $x_0-\sum_{j<k}\sigma_jr_j(1,1)$, $k\ge0$, obtained by feeding the rule a finite sequence $(g_0,\dots,g_{k-1})$ of these eight vectors, where $\sigma_j$ is the common sign of the components of $g_j$ and $r_j$ is the radius the rule prescribes after observing $g_0,\dots,g_j$. There are at most $8^k$ such sequences of length $k$, so $P$ is countable, and $P$ does not depend on $c$. Let $N$ be the set of all $c$ such that $c_1=p_1$, or $c_2=p_2$, or $c_1+c_2=p_1+p_2$ for some $p\in P$, or $c_1-c_2=x_{0,1}-x_{0,2}$. Then $N$ is a countable union of lines and has Lebesgue measure zero.

Let $c\notin N$. Induction on $k$ shows that every iterate $x_k$ belongs to $P$ and is a point of differentiability of $f_c$: this holds for $x_0\in P$, and if it holds for $x_0,\dots,x_k$, then the observed subgradients are the gradients $\nabla f_c(x_0),\dots,\nabla f_c(x_k)$ whatever the selection, so $x_{k+1}\in P$, and $x_{k+1}$ lies on none of the three lines because $c\notin N$. Hence the iterates are the same for every subgradient selection, and $x_{k,1}-x_{k,2}=x_{0,1}-x_{0,2}$ for all $k$ because every step is a multiple of $(1,1)$. Finally, $\min f_c=0$ and $f_c(x)\ge\varepsilon(|x_1-c_1|+|x_2-c_2|)\ge\varepsilon|(x_1-c_1)-(x_2-c_2)|$, which gives the lower bound of Remark~\ref{rem:ce}; it is positive because $c_1-c_2\ne x_{0,1}-x_{0,2}$.

\bibliographystyle{siamplain}
\bibliography{references}

\end{document}